\documentclass[aap]{imsart}

\usepackage[square,compress,comma, numbers,sort]{natbib}
\usepackage[colorlinks=true, citecolor=red, linkcolor=blue]{hyperref}
\usepackage{amsfonts,mathtools}

\usepackage{amsmath}
\usepackage{bbm}
\usepackage{amssymb}
\usepackage{mathtools}

\numberwithin{equation}{section}

\usepackage{color}

\DeclareMathOperator{\one}{\mathbbm{1}}

\definecolor{c20}{rgb}{0.,0.7,0.}
\definecolor{c30}{rgb}{0.,0.,1.}
\definecolor{c40}{rgb}{1,0.1,0.7}
\definecolor{c50}{rgb}{1,0,0}
\definecolor{c60}{rgb}{1,0.9,0.1}

\def\hfr#1{\mathfrak{h}_{#1}}
\def\mfr#1{\mathfrak{m}^{#1}}

\newcommand{\abs}[1]{\left\lvert #1 \right\rvert}

\newcommand{\sprod}[1]{\left\langle#1\right\rangle}

\newcommand{\E}[1]{\mathbb{E}\left\{ #1\right\}}

\newcommand{\pk}[1]{\mathbb{P} \left\{ #1 \right \} }

\newcommand{\R}{\mathbb{R}}

\newcommand{\N}{\mathbb{N}}

\newcommand{\BQN}{\begin{eqnarray}}
\newcommand{\EQN}{\end{eqnarray}}
\newcommand{\BQNY}{\begin{eqnarray*}}
\newcommand{\EQNY}{\end{eqnarray*}}

\newcommand{\BS}{\begin{sat}}
\newcommand{\ES}{\end{sat}}
\newcommand{\BT}{\begin{theo}}
\newcommand{\ET}{\end{theo}}
\newcommand{\BK}{\begin{korr}}
\newcommand{\EK}{\end{korr}}

\newcommand{\BD}{\begin{de}}
\newcommand{\ED}{\end{de}}

\newcommand{\BIT}{\begin{itemize}}
\newcommand{\EIT}{\end{itemize}}
\newcommand{\BDI}{\begin{description}}
\newcommand{\EDI}{\end{description}}

\newcommand{\BRM}{\begin{remarks}}
\newcommand{\ERM}{\end{remarks}}

\newcommand{\BEL}{\begin{lem}}
\newcommand{\EEL}{\end{lem}}

\newtheorem{theo}{Theorem}[section]
\newtheorem{sat}[theo]{Proposition}
\newtheorem{de}[theo]{Definition}
\newtheorem{lem}[theo]{Lemma}

\newtheorem{korr}[theo]{Corollary}
\newtheorem{remark}[theo]{Remark}
\newtheorem{remarks}[theo]{Remarks}
\newtheorem{prop}[theo]{Proposition}

\newcommand{\nelem}[1]{{Lemma \ref{#1}}}

\newcommand{\neprop}[1]{{Proposition \ref{#1}}}
\newcommand{\netheo}[1]{{Theorem \ref{#1}}}

\newcommand{\prooftheo}[1]{\noindent\textsc{\bf Proof of Theorem} \ref{#1}:}
\newcommand{\proofprop}[1]{\noindent\textsc{\bf Proof of Proposition} \ref{#1}:}
\newcommand{\prooflem}[1]{\noindent\textsc{\bf Proof of Lemma} \ref{#1}:}

\newcommand{\COM}[1]{}

\newcommand{\CN}[1]{{\color{red}#1}}

\def\td{\text{\rm d}}

\newcommand{\QED}{\hfill $\Box$}

\newcommand{\kb}[1]{\boldsymbol{#1}}
\newcommand{\vk}[1]{\kb{#1}}
\newcommand{\LON}{\lambda_{\mathrm{on}}}
\newcommand{\LOFF}{\lambda_{\mathrm{off}}}

\begin{document}
	\begin{frontmatter}
\title{Functional Central Limit Theorem for the principal eigenvalue of dynamic Erd\H{o}s-R\'enyi random graphs}
\runtitle{Principal eigenvalue of dynamic ERRG}
\begin{aug}
\author[A]{\fnms{Rajat Subhra} ~\snm{Hazra}\ead[label=e1]{r.s.hazra@math.leidenuniv.nl}}
\author[B]{\fnms{Nikolai}~\snm{ Kriukov}\ead[label=e2]{n.kriukov@uva.nl}}
\author[A, B]{\fnms{Michel}~\snm{ Mandjes}\ead[label=e3]{m.r.h.mandjes@math.leidenuniv.nl}}
\address[A]{Mathematical Institute, Leiden University, P.O. Box 9512,
2300 RA Leiden,
The Netherlands\printead[presep={,\ }]{e1,e3}}
\address[B]{Korteweg-de Vries Institute for Mathematics, University of Amsterdam, Science Park 904, 1098 XH Amsterdam, The Netherlands\printead[presep={,\ }]{e2,e3}}
\end{aug}
\date{\today}

\maketitle

\begin{abstract}
    In this paper we consider a dynamic version of the Erd\H{o}s-R\'{e}nyi random graph, in which edges independently appear and disappear in time, with the on- and off times being exponentially distributed. 
    The focus lies on the evolution of the principle eigenvalue of the adjacency matrix in the regime that the number of vertices grows large. The main result is a functional central limit theorem, which displays that the principal eigenvalue essentially inherits the characteristics of the dynamics of the individual edges.
\end{abstract}


\noindent
%
%
%
%

\end{frontmatter}

\section{Introduction}

    Over the past few decades a broad range of {\it random graph models} has been proposed and analyzed. These studies aim to capture the essentials of various types of real-life networks, inspired by applications in e.g.\  sociology, biology, economics, physics, climate, computer science, and epidemiology. Arguably the most intensively studied instance is the \textit{Erd\H{o}s-R\'{e}nyi random graph} \cite{erdos1959random, gilbert1959random}. This model is characterized through a fixed amount of vertices $N$, where any possible vertex pair is connected with a given fixed probability~$p$ (independently from every other vertex pair).
    
    When considering the random graph literature, it is observed that most of it concerns \textit{static} models, in which the graph is sampled once and does not evolve thereafter. However, virtually any real-life network has a structure that changes over time: typically, edges may appear or disappear. This has motivated a growing interest in \textit{dynamic} random graph models, and it in particular gave rise to the development of a  dynamic version of the Erd\H{o}s-R\'{e}nyi random graph; see e.g.\ \cite{zhang2017random, mandjes2019dynamic, holme2012temporal, holme2015modern, braunsteins2023sample}. In the dynamic random graph literature, the main objective lies in the analysis of various time-dependent properties.

    \medskip

    Perhaps the most elementary dynamic Erd\H{o}s-R\'{e}nyi random graph model, studied in e.g.\ \cite{zhang2017random,mandjes2019dynamic, braunsteins2023sample}, is defined as follows. The initial graph, pertaining to time $0$, is sampled as a static Erd\H{o}s-R\'{e}nyi random graph on $N$ vertices with edges independently sampled with probability~$p$. After that, edges independently alternate between being absent and present, with the `on-times' being exponentially distributed with parameter $\LON$,  and the `off-times' exponentially distributed with parameter $\LOFF$. This dynamic random graph mechanism was considered in \cite{zhang2017random}, but this paper had a broader focus: the objective was to propose dynamic random graph models whose stationary behavior reproduces that of a series of known static random graph models, including related parameter inference issues. The main goal of \cite{mandjes2019dynamic} was to introduce a family of dynamic Erd\H{o}s-R\'{e}nyi random graph models (of which the above `exponential on-off version' is a special case) and to prove that in a specific scaling regime the number of edges (as a function of time, that is) converges to an Ornstein-Uhlenbeck process. Reference \cite{braunsteins2023sample} establishes the sample-path large deviations of the exponential on-off Erd\H{o}s-R\'{e}nyi random graph, thus extending the static counterpart developed in the seminal work \cite{chatterjeevaradhan2011static}. In principle any quantity or phenomenon that has been analyzed in the static random graph literature (for instance properties of a giant component, or the number of triangles or any other subgraph count) can now be studied for its dynamic counterpart.

\medskip

     As pointed out in e.g.\ \cite{castellano2017relating, pastor2018eigenvector, martin2014localization, newman2006finding}, in order to get a handle on for instance the network's centrality and on its community structures, insight into {\it spectral properties} of the random graph's adjacency matrix are of great importance. Motivated by this observation, in \cite{10.1214/11-AOP734} a central limit theorem (CLT) has been established for the principal eigenvalue in the static Erd\H{o}s-R\'{e}nyi model; see also recent papers \cite{dionigi2023central, chakrabarty2:hazra} which provide the analogous results for the Chung-Lu model and inhomogeneous Erd\H{o}s-R\'enyi random graphs, respectively. These results concern the regime in which the number of vertices grows large: after an appropriate centering and normalization, the principal eigenvalue tends to a Gaussian random variable as $N\to\infty$. 
     The main goal of the present paper is to extend these existing static results to the dynamic case: we aim to analyze the time-dependent behavior of the principal eigenvalue for the exponential on-off Erd\H{o}s-R\'{e}nyi dynamic random graph model that we introduced above. Our concrete objective is to establish a {\it functional} central limit theorem (FCLT) in the regime that $N$ grows large.

 \medskip

The main principle underlying the analysis for the static Erdős-Rényi model in \cite{10.1214/11-AOP734}, is that the principal eigenvalue essentially inherits the probabilistic properties of the total number of edges; evidently, this total number of edges behaves (after an appropriate centering and scaling) as a Gaussian random variable as $N \to \infty$ by virtue of the conventional CLT. In the present paper, we extend this principle to the setting of dynamic Erd\H{o}s-R\'enyi models. Our result reveals that the principal eigenvalue process, in the limiting regime considered, converges to a Gaussian process that has the same correlation structure as the edge processes.

At the methodological level, in the study of the principal eigenvalue of the adjacency matrix a key role is played by rank-one perturbations of Wigner-type matrices. There is a vast body of literature concerning finite-rank perturbations of Wigner matrices, which have proven to be useful in various statistical applications such as community detection or the detection of denser communities in random graphs \cite{AriasCastro, Abbe, Fan}. A more extensive account of the literature is provided in Section \ref{section:main results}.

To the best of our knowledge, there has been no study on the {\it time-dependent} behavior of the eigenvalues of {\it dynamic} random graphs. In our first main result, Theorem \ref{main}, we state the FCLT for the principal eigenvalue of our dynamic Erdős-Rényi model. This result crucially relies on establishing a representation of the principal eigenvalue uniformly over time, which is achieved in Theorem \ref{true_eigenvalue_representation}. Our line of reasoning essentially follows the one used in \cite{10.1214/11-AOP734}, but a series of specific challenges need to be resolved, typically related to the fact that we need that certain properties relied upon in \cite{10.1214/11-AOP734} extend to the time-dependent setting. Using the representation of Theorem~\ref{true_eigenvalue_representation}, we succeed in proving the FCLT by first establishing the finite-dimensional convergence and then proving tightness in the space of càdlàg functions. In our work we do not let the model parameters depend on $N$; see Remark \ref{R1a} for a short account of the complications that would arise when making these parameters $N$-dependent.

    \medskip
    
    This paper is organized as follows. In Section 2 we formally describe the model of our exponential on-off Erd\H{o}s-R\'{e}nyi dynamic graph, also presenting a number of basic properties of the edge processes. In Section 3 we state the main results of our paper: \netheo{main} presents the FCLT of the principal eigenvalue to a specific Gaussian process, \netheo{true_eigenvalue_representation} presents a uniform `with high probability representation' for this eigenvalue, and \netheo{thm_expectation} gives the uniform asymptotics for the expectation of the principal eigenvalue. Section~4 covers the proof of \netheo{main}, Section~5 the proof of \netheo{true_eigenvalue_representation}, and  Section~6 the proof of \netheo{thm_expectation}. Section 7 provides concluding remarks and discusses suggestions for follow-up research.
    
\subsection*{Acknowledgments} 
This research was supported by the European Union’s Horizon 2020 research and innovation programme under the Marie Sklodowska-Curie grant agreement no.\ 945045, and by the NWO Gravitation project NETWORKS under grant agreement no.\ 024.002.003. \includegraphics[height=1em]{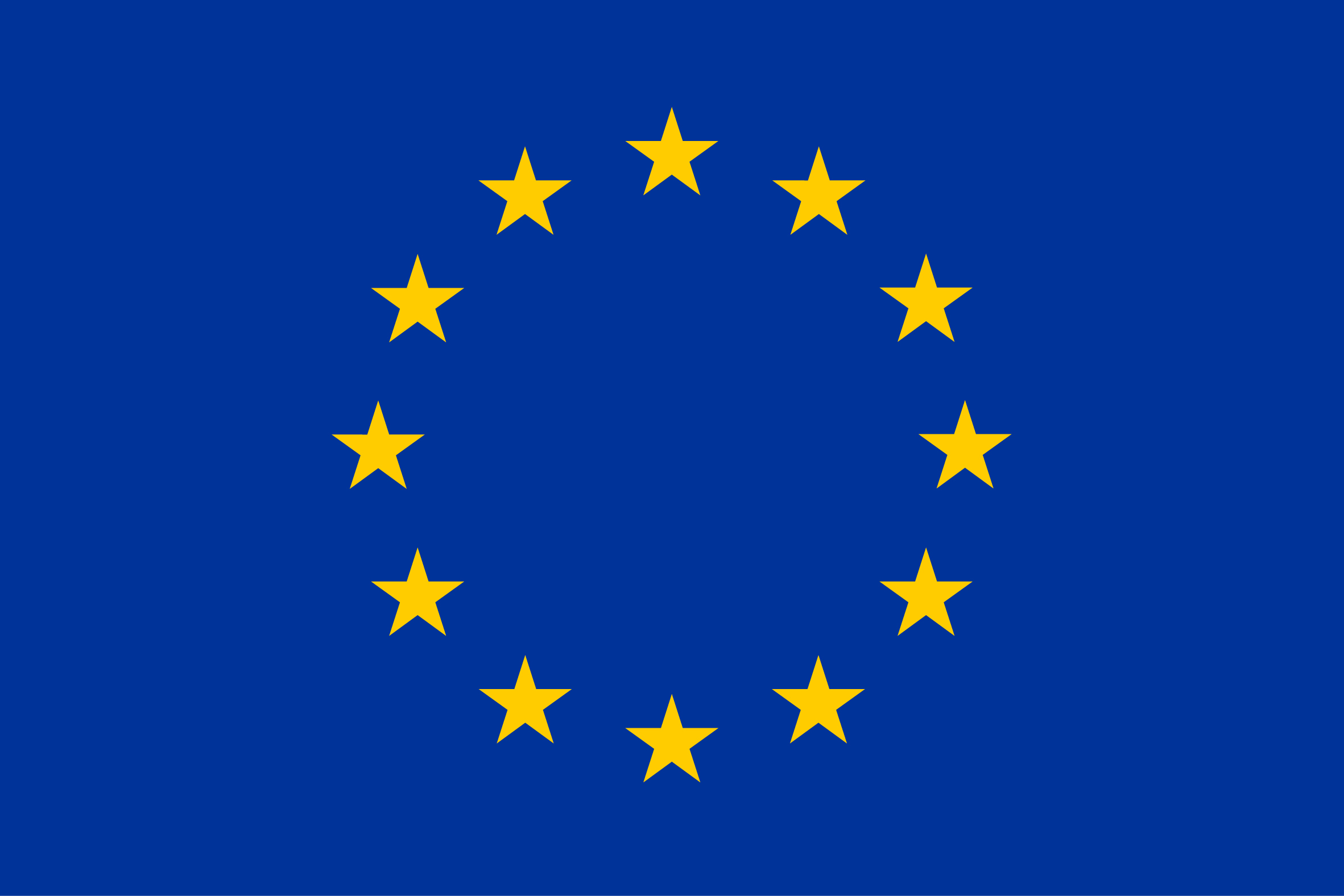}

\section{Model and preliminaries}

Consider a {\it dynamical Erd\H{o}s-R\'{e}nyi random graph} $G_N(t)$ with $N$ vertices over the time interval $t\in[0,T]$ for some $T>0$, evolving as follows. 
In the initial graph ${G}_N(0)$ for each vertex pair $(i,j)$ with probability $p(0)\in(0,1)$ an edge is present, where each of these potential edges is sampled independently.
After time 0, each edge evolves independently according to the following rules: if some edge exists at time $t\geqslant 0$, it remains present for a exponentially distributed time with parameter $\LON>0$, and then disappears. 
Otherwise, if at time $t\geqslant 0$ some edge does not exist, it is absent for an exponentially distributed time with parameter $\LOFF>0$, after which it appears. 

Fixing a time horizon $T>0$, the adjacency matrix process $(A_N(t))_{t\in[0,T]}$ of this graph process $({G}_N(t))_{t\in[0,T]}$ can be formally represented through the {$N^2$} processes $(a_{i,j}(t))_{t\in[0,T]}$, defined by means of the following three rules. Throughout this paper we use the definition $[n]:=\{1,\ldots,n\}$.
\begin{itemize}
\item[$\circ$]
At time $t=0$ the edges of the graph are initialized via 
\begin{align*}
a_{i,j}(0)&=\begin{cases}
    1\quad\text{with probability } p(0),\\
    0\quad\text{with probability } \bar p(0):=1-p(0),
\end{cases}
\end{align*}
with $j\in[N]$ and $i\in[j]$, where each of these $a_{i,j}(0)$ is sampled independently. 
\item[$\circ$]
After time $0$, the processes $(a_{i,j}(t))_{t\in[0,T]}$, with $j\in[N]$ and $i\in[j]$, evolve independently, as follows. Fix $i,j$ such that 
$j\in[N]$ and $i\in[j]$.
Let $\bigl(\xi^1_{\ell}(i,j)\bigr)_{\ell \in {\mathbb N}}$ be a sequence of iid  (independent and identically distributed) exponentially distributed random variables with parameter $\LON>0$, and analogously $\bigl(\xi^0_{\ell}(i,j)\bigr)_{\ell \in {\mathbb N}}$ a sequence of iid exponentially distributed random variables with parameter $\LOFF>0$. Define the duration of the first $k$ on-off cycles (for $k\in{\mathbb N}$) by
\[T_k(i,j)\equiv T_k:=\sum_{\ell = 1}^k \xi^0_\ell(i,j)+\sum_{\ell = 1}^k \xi^1_\ell(i,j),\]
with $T_0(i,j)\equiv T_0:=0$.
If $a_{i,j}(0)=0$, then $a_{i,j}(t)=0$ if 
\[t\in[T_k,T_k+\xi^0_{k+1}(i,j))\] for some $k\in {\mathbb N}_0$, and $a_{i,j}(t)=1$ otherwise.
On the other hand, if $a_{i,j}(0)=1$, then $a_{i,j}(t)=1$ if 
\[t\in[T_k,T_k+\xi^1_{k+1}(i,j))\] for some $k\in {\mathbb N}_0$, and $a_{i,j}(t)=0$ otherwise. 
\item[$\circ$]
Above we defined $(a_{i,j}(t))_{t\in[0,T]}$ for $j\in[N]$ and $i\in[j]$.
As a last step, we define the remaining processes $(a_{i,j}(t))_{t\in[0,T]}$. The evolving graph is non directed, entailing that at any point in time the adjacency matrix is symmetric. This means that $a_{i,j}(t) = a_{j,i}(t)$, for $j\in[N]$ and $i\in[N]\setminus [j]$.
\end{itemize}
For later use, we also introduce the following two processes: for $t\geqslant 0$,
\begin{align}
    \begin{split}
    \mathcal{P}_{0,i,j}(t) &:=\max\left\{k\in\N_0 \colon \sum_{\ell=1}^{\lceil\frac{k}{2}\rceil} \xi^0_{\ell}(i,j) + \sum_{\ell=1}^{\lfloor\frac{k}{2}\rfloor} \xi^1_{\ell}(i,j)\leqslant t\right\},\\
    \mathcal{P}_{1,i,j}(t) &:=1 + \max\left\{k\in\N_0 \colon \sum_{\ell=1}^{\lfloor\frac{k}{2}\rfloor} \xi^0_{\ell}(i,j) + \sum_{\ell=1}^{\lceil\frac{k}{2}\rceil} \xi^1_{\ell}(i,j)\leqslant t\right\}.
    \end{split}\label{P_def}
\end{align}
The random variables $\xi^0_{\ell}$ and $\xi^1_{\ell}$ for $\ell=0$ are set to be $0$.

It can be verified that each of the processes $(a_{i,j}(t))_{t\in[0,T]}$ can be alternatively represented as
\[a_{i,j}(t)=\frac{1}{2}-\frac{1}{2}(-1)^{{\mathcal P}_{a_{i,j}(0),i,j}(t)}.\]
To see this, note that (for instance) $a_{i,j}(0)=0$ and $\mathcal{P}_{0,i,j}(t)$ being odd means that there was no edge at time $0$ but that there is an edge at time $t$; the other three cases can be dealt with analogously. 

\begin{remark}{\em 
    Notice that in our setup, {\it self-loops} are allowed: $(a_{i,i}(t))_{t\in[0,T]}$ is distributionally identical to $(a_{i,j}(t))_{t\in[0,T]}$ for $i\not = j$. This entails that, effectively, $A_N(t)$ is described by {$\frac{1}{2}N(N+1)$} of its entries. We could have alternatively studied the model {\it without} self-loops, with only {$\frac{1}{2}N(N-1)$} relevant entries. By inspecting all proofs, the reader can verify that our results carry over to the counterpart of our model that has no self-loops.} $\hfill\Diamond$
\end{remark}

\medskip

We proceed by presenting, for later reference, a number of elementary properties of the per-edge processes.
To this end, we consider the probability, for $t\in[0,T]$,
\begin{align}
\label{p}
    p(t) := \pk{a_{i,j}(t) = 1}= \bar p(0)\,p_{01}(t) + p(0)\,p_{11}(t),
\end{align}
where, for $k,\ell\in\{0,1\}$, we define $p_{k\ell}(t) := \pk{a_{i,j}(t)=\ell\,|\, a_{i,j}(0)=k}.$
It follows by standard techniques that, with $\varrho:=\LOFF/(\LON+\LOFF )$,
\begin{align*}
    p_{01}(t) = \varrho  - \varrho \,e^{-(\LON+\LOFF )t},\:\:\:\:
    p_{11}(t) = \varrho + (1-\varrho)\,e^{-(\LON+\LOFF )t},
\end{align*}
so that, after a straightforward computation in which we use \eqref{p},
\begin{align}
    \label{p_def} p(t) = \varrho + \big((1-\varrho) \,p(0) - \varrho\,\bar p(0)\big)\,e^{-(\LON+\LOFF )t}.
\end{align}
In particular, the processes $(a_{i,j}(t))_{t\in[0,T]}$ are continuous-time Markov processes.
A property which we heavily rely upon in this paper is that, as a consequence of $p(0)\in(0,1)$, there exist constants $p^-,p^+\in(0,1)$ such that
\begin{align}
p^-\leqslant \inf_{t\in[0,T]}p(t) \leqslant \sup_{t\in[0,T]}p(t) \leqslant p^+.\label{p_bounds}
\end{align}
For any $i,j$, using \eqref{p_def} and the Markov property, the per-edge covariance can be determined. Indeed, for $0\leqslant t_1\leqslant t_2\leqslant T$,
    \begin{align*}
        \operatorname{Cov}(a_{i,j}(t_1),a_{i,j}(t_2))&=\E{(a_{i,j}(t_1)-p(t_1))(a_{i,j}(t_2)-p(t_2))}\\
        &=\E{a_{i,j}(t_1)\,a_{i,j}(t_2)}-p(t_1)\,p(t_2)\\
        &=\pk{a_{i,j}(t_1)=1,\,a_{i,j}(t_2)=1}-p(t_1)\,p(t_2)\\
        &=p(t_1)\big(p_{11}(t_2-t_1) - p(t_2)\big).
    \end{align*}
Then observing that $p(t_2)=p(t_1)\,p_{11}(t_2-t_1) + (1-p(t_1))\,p_{01}(t_2-t_1)$ (again by the Markov property), we thus obtain
    \begin{align}
        \operatorname{Cov}(a_{i,j}(t_1),a_{i,j}(t_2))
        &=p(t_1)\bigl(p_{11}(t_2-t_1) -p(t_1)p_{11}(t_2-t_1) - (1-p(t_1))p_{01}(t_2-t_1)\bigr)\notag\\
        &=p(t_1)\bigl((1-p(t_1))p_{11}(t_2-t_1) - (1-p(t_1))p_{01}(t_2-t_1)\bigr)\notag\\
        &=p(t_1)(1-p(t_1))\left(p_{11}(t_2-t_1) - p_{01}(t_2-t_1)\right)\notag\\
        &=p(t_1)(1-p(t_1))\,e^{-(\LON+\LOFF )(t_2-t_1)}.\label{eq:COV}
    \end{align}

\section{Main results}\label{section:main results}

As mentioned in the introduction, the primary objective of this paper lies in characterizing the time-dependent behavior of the principal eigenvalue $\mu_{N}(t)$ of the adjacency matrix $A_N(t)$, leading to an FCLT for $(\mu_N(t))_{t\in[0,T]}$. In this section we present our three main theorems, but we start by providing some background, thus putting our findings into perspective.

Quantitative estimates of the principal eigenvalue of the homogeneous, \emph{static} Erd\H{o}s-R\'enyi random graph on $N$ vertices with parameter $p_N$  (when $Np_N\gg (\log N)^4$) have been identified in \cite{furedi1981eigenvalues,vu2007spec}. It followed from these works that the smallest and second-largest eigenvalues converge to the edge of the support of the semicircular law. These results were further improved by  \cite{benaych2020spectral}, where the condition on sparsity was extended to the case $Np_N\gg \log N$ (which is also the connectivity threshold); it was in addition shown that {\it inhomogeneous} static random graphs exhibit a similar behavior. Under the assumption that $N^{\xi}\ll Np_N$ for some $\xi\in (2/3, 1]$, it was shown in \cite[Theorem 2.7]{Erdos2} that the second-largest eigenvalue of the (homogeneous, static) Erd\H{o}s-R\'enyi random graph, after an appropriate centering and scaling, converges in distribution to the Tracy-Widom law (results that were recently improved on in \cite{lee2018}). Properties of the principal eigenvector (i.e., the eigenvector corresponding to the principal eigenvalue) in the homogeneous case were studied in \cite{tran2013sparse, Erdos2, lee2018}.

A result that is of particular importance in the context of the present paper is the following. In \cite{10.1214/11-AOP734} it was shown that the principal eigenvalue of the static Erd\H{o}s-R\'enyi random graph has a Gaussian limit. Indeed, as $N\to\infty$, 
\begin{equation}\label{univ}(2\,p_N(1-p_N))^{-1/2}(\mu_N- \E{\mu_N})\overset{\rm d}\rightarrow X,\end{equation}
with $X$ a standard normal random variable. The above result was shown under the assumption that $(\log N)^{\xi} \ll Np_N$ for some $\xi>8$ (entailing that the result holds in the specific case that $p$ does not depend on $N$).

It is well known that in the classical case of a (standard) Wigner matrix, the principal eigenvalue converges to the Tracy-Widom law. We note that there is a different scaling between the edge and bulk of the spectrum in Erd\H{o}s-R\'enyi random graphs: the bulk scales at rate $(Np_N(1-p_N))^{1/2}$, whereas the principal eigenvalue has the scaling $Np_N(1-p_N)$. Letting \begin{equation} \label{eq.defw} H_N:=A_N-\E{A_N}, \end{equation} where $\E{A_N}$ is the entrywise expectation of $A_N$, it can easily be seen that
\[A_N= p_N \one\one^\prime + H_N,\]
with $\one$ denoting the $N\times1$ vector with each entry equalling $1$. In the static case, the empirical spectral distribution of $(Np_N(1-p_N))^{-1/2}H_N$ converges to the semicircle law, the principal eigenvalue of the same converges to $2$ almost surely.
As $\E{A_N}$ is a rank-one matrix, it turns out that the principal eigenvalue of $A_N$ scales like $Np_N$, which is different from the bulk scaling.

To analyze the fluctuations of the principal eigenvalue of $A_N$ around its mean, one needs to study in detail the effects of rank-one perturbations on a Wigner-type matrix $H_N$. When $H_N$ is a symmetric random matrix with independent and identically distributed entries and the perturbation comes from a rank-one matrix, then the fluctuation of the principal eigenvalue depends on the form of the {\it deformation matrix} \cite{feral2007largest, capitaine2009largest, capitaine2012central}.
 The rank-one case was extended to the finite-rank case in the works  \cite{benaych:2011, pizzo:2013, chakrabarty2:hazra}.

 The adjacency matrix of our dynamic random graph mechanism does not fall directly into the purview of the above results. It is true that $H_N$, as given by \eqref{eq.defw}, is a symmetric matrix, with independent entries above the diagonal, but the crucial issue is that in our setup the entries are now functions of time. It is this time dependence that prevents us from the (direct) use of results and estimates that were developed earlier for the static case.
 
 Our main result concerns the following FCLT for $(\mu_{N}(t))_{t\in[0,T]}$. Specializing to our framework in which the per-edge processes do not depend on $N$, Theorem \ref{main} generalizes the univariate CLT for the static Erd\H{o}s-R\'enyi random graph that was given above, i.e., claim (6.12) in \cite[Theorem 6.2]{10.1214/11-AOP734}; compare the claim of the below theorem with the univariate version \eqref{univ}.

\begin{theo}\label{main}
    As $N\to\infty$
    \begin{align*}
        \mu_{N}(t)-\E{\mu_{N}(t)}\xrightarrow[]{\rm d}  X(t)
    \end{align*}
    in $\mathbb{D}([0,T],\R)$, where $\bigl(X(t)\bigr)_{t\in[0,T]}$ is a centered Gaussian process characterized by the covariance function, for $0\leqslant t_1\leqslant t_2\leqslant T$,
    \begin{align*}
        \operatorname{Cov}(X(t_1),X(t_2))=2\,p(t_1)(1-p(t_1))\,e^{-(\LON+\LOFF )(t_2-t_1)}.
    \end{align*}
\end{theo}


Our proof borrows elements from the proof of claim (6.12) in \cite[Theorem 6.2]{10.1214/11-AOP734}, but along the way various technicalities need to be resolved. In particular the first step towards showing the FCLT is to establish a uniform version of \cite[Eqn.\ (6.11)]{{10.1214/11-AOP734}}, so as to deal with the time-dependent nature that is inherent to our dynamic setup; a crucial role is played by a representation that is derived from the series expansion \eqref{mu_recurrent_expansion}. Importantly, even though the results presented in our work (see in particular our \netheo{true_eigenvalue_representation}) may look similar to their counterparts in \cite{{10.1214/11-AOP734}}, one cannot derive them straightforwardly from results in \cite{{10.1214/11-AOP734}}. The main difficulty lies in the fact that in \cite{10.1214/11-AOP734} the representation for the principal eigenvalue is not given in an `almost-sure' sense, but just `with high probability'. This in particular entails that such a representation is satisfied on an event of which the probability is not equal to one, but tends to one sufficiently fast as $N$ grows. As one cannot derive an almost sure representation for the principal eigenvalue, we need to establish one that is uniform in two ways: {not only need the remainder terms in the series expansion \eqref{mu_recurrent_expansion} be uniformly bounded}, but also the high probability event should be uniform over all points of time. After having settled this uniform representation, we need to go through the usual steps to prove the FCLT: first we are to prove convergence of the finite-dimensional distributions to those of the stated Gaussian limit, and then we have to verify tightness conditions (which poses a number of technical challenges).

\begin{theo}\label{true_eigenvalue_representation}
    There exist a set of events $\Omega_N$ such that for some $\eta>0$
    \begin{align*}
        \pk{\Omega^{\rm c}_{N}}\leqslant e^{-\eta (\log N)^2},
    \end{align*}
    while on $\Omega_N$ for all $t\in[0,T]$
    \begin{align}\label{eq:DEC}
        \mu_{N}(t)=\E{\mu_{N}(t)} + \frac{1}{N}\sum_{i,j=1}^{N}\left(a_{i,j}(t)-p(t)\right) + \mathcal{G}_N(t)\left(\frac{(\log N)^{4}}{\sqrt{N}}\right),
    \end{align}
    for some ${\mathcal G}_N(t)$ such that $\limsup_{N\to\infty}\sup_{t\in[0,T]}\abs{\mathcal{G}_N(t)}$ is finite almost surely.
\end{theo}

\begin{remark}\label{R1}
    {\em Observe the similarity between the covariance structure of $(X(t))_{t\in[0,T]}$, and the covariance pertaining to a single edge as given by \eqref{eq:COV}. One could say that the principal eigenvalue process essentially inherits the correlation structure of the edges (where the factor $2$ is a consequence of the symmetry of the adjacency matrix), which could have been anticipated in view of \eqref{eq:DEC}. Observe also that
    \[\lim_{t\to\infty} \operatorname{Cov}(X(t),X(t+\delta))=2\,\varrho(1-\varrho)\,e^{-(\LON+\LOFF )\delta}, \]
    where we recognize the correlation structure of the (stationary version of the) Ornstein-Uhlenbeck process. 
    \hfill$\Diamond$}
\end{remark}

In Theorem \ref{main}, we did not explicitly assess the centering of the principal eigenvalue. The next theorem gives an asymptotic expansion of the expectation $\E{\mu_N(t)}$ in the regime $N\to\infty.$
\begin{theo}\label{thm_expectation} The expectation of the principal eigenvalue $\mu_N(t)$ satisfies, for $t\in[0,T]$,
    \begin{align*} 
        \E{\mu_N(t)} = Np(t)+(1-p(t)) + \frac{\overline{\mathcal{G}}_N(t)}{N},
    \end{align*}
    for some $\overline{\mathcal G}_N(t)$ such that  $\limsup_{N\to\infty}\sup_{t\in[0,T]}\abs{\overline{\mathcal{G}}_N(t)}$ is finite.
\end{theo}

\begin{remark}\label{R1a}{\em 
In our work we consider the case of $p_N(0)$ being independent of $N$. To let $p_N(0)$ depend on $N$, a natural approach would be to work with $N$-dependent edge processes characterized by transition rates $\lambda_{{\rm on},N}$ and $\lambda_{{\rm off},N}$. While the expression for the covariance, as given in \eqref{eq:COV}, remains valid, for the covariance function to have a non-trivial limit, the edge processes require time-rescaling. This introduces multiple hurdles in proving tightness, as well as in the proof of Proposition~\ref{model_basic_properties}. To align with the estimates used in \cite{{10.1214/11-AOP734}} and this article, we would need to impose the restriction that $Np_N(0)$ is much larger than $(\log N)^\xi$, for some $\xi\geqslant 1$. This extension constituting various additional technical challenges, we do not explore it in this article. 
    \hfill$\Diamond$}
\end{remark}

\section{Proof of Theorem \ref{main}}

In this section we prove Theorem \ref{main}, assuming Theorem~\ref{true_eigenvalue_representation} to hold (which we prove in the next section). We remark that, compared to the results for the static case presented in \cite{10.1214/11-AOP734}, in our dynamic random graph setting there are the complications inherent to {\it functional} convergence; most notably, as pointed out below, for Theorem \ref{main} to hold, tightness has to be validated. 

\medskip

\prooftheo{main} 
  To show the process-level distributional convergence, we follow the standard procedure of first showing the convergence of finite-dimensional distributions and then establishing tightness. For the first point, fix some $0\leqslant t_1<\ldots<t_m\leqslant T$ for some $m\in{\mathbb N}$. To show that vector $$(\mu_{N}(t_1)-\E{\mu_{N}(t_1)},\ldots,\mu_{N}(t_m)-\E{\mu_{N}(t_m)})$$ converges to a Gaussian vector we use the Cram\'er-Wold device: it is enough to show that, for any vector ${\boldsymbol q}\equiv (q_1,\ldots,q_m)\in\R^m$, the random variable $\sum_{i=1}^m q_i(\mu_{N}(t_i)-\E{\mu_{N}(t_i)})$ converges to a Normally distributed random variable. Using \netheo{true_eigenvalue_representation}, there is a set of events $\Omega_N$ such that
    \begin{align*}
        \pk{\Omega^{\rm c}_N}\leqslant e^{-\eta\, (\log N)^2}, 
    \end{align*}
    where on $\Omega_N$, for ${\mathcal G}_N(t)$ as specified in the theorem,
    \begin{align*}
        &\sum_{i=1}^m q_i(\mu_{N}(t_i)\,-\,\E{\mu_{N}(t_i)})\\
        &\qquad=\frac{1}{N}\sum_{i=1}^{m}q_i\sum_{j,k=1}^{N}\left(a_{j,k}(t_i)-p(t_i)\right) + \sum_{i=1}^{m}q_i\,\mathcal{G}_N(t_i)\left(\frac{(\log N)^4}{\sqrt{N}}\right)\\
        &\qquad=\frac{1}{N}\sum_{j,k=1}^{N}\left(\sum_{i=1}^{m}q_i(a_{j,k}(t_i)-p(t_i))\right)+\sum_{i=1}^{m}q_i\,\mathcal{G}_N(t_i)\left(\frac{(\log N)^4}{\sqrt{N}}\right).
    \end{align*}
    Using that $a_{j,k}(t_i) = a_{k,j}(t_i)$,  the fact that the random variables  $\sum_{i=1}^{m}q_i(a_{j,k}(t_i)-p(t_i))$ are independent for different pairs $(j,k)$ with $j<k$, and $\limsup_{N\to\infty}\sup_{t\in[0,T]}\abs{\mathcal{G}_N(t)}<\infty$ (by virtue of \netheo{true_eigenvalue_representation}), 
    \begin{align*}
        &\lim_{N\to\infty}\pk{\frac{1}{N}\sum_{j,k=1}^{N}\left(\sum_{i=1}^{m}q_i(a_{j,k}(t_i)-p(t_i))\right)+\sum_{i=1}^{m}q_i\,\mathcal{G}_N(t_i)\left(\frac{(\log N)^4}{\sqrt{N}}\right)<x}\\
        &=\lim_{N\to\infty}\pk{\frac{1}{N}\sum_{j,k=1}^{N}\left(\sum_{i=1}^{m}q_i(a_{j,k}(t_i)-p(t_i))\right)<x}\\
        &=\lim_{N\to\infty}\pk{\frac{2}{N}\sum_{\substack{j,k=1\\j<k}}^{N}\left(\sum_{i=1}^{m}q_i(a_{j,k}(t_i)-p(t_i))\right) + \frac{1}{N}\sum_{j=1}^{N}\sum_{i=1}^m q_i(a_{jj}(t_i)-p(t_i))<x}\\
        &=\pk{\xi<x},
    \end{align*}
    with $\xi$ denoting a centered normal random variable with variance $2\sum_{i=1}^{m}q_i^2\operatorname{Var}\bigl(a_{1,1}(t_i)\bigr)$, where in the last step the conventional central limit theorem has been used (also observing that the sum over $j<k$ consists of $N(N-1)/2$ terms).
    Then using that, by \netheo{true_eigenvalue_representation},
    \begin{align*}
    &\left\lvert\:\pk{\sum_{i=1}^m q_i(\mu_{N}(t_i)-\E{\mu_{N}(t_i)})<x}\right.\\
    &\qquad\qquad\left.-\:\pk{\frac{1}{N}\sum_{j,k=1}^{N}\left(\sum_{i=1}^{m}q_i(a_{j,k}(t_i)-p(t_i))\right)+\sum_{i=1}^{m}q_i\,\mathcal{G}_N(t_i)\left(\frac{(\log N)^4}{\sqrt{N}}\right)<x}\,\right\rvert\\
    &\qquad\qquad\leqslant 2\,\pk{\Omega^{\rm c}(N)}\leqslant 2\,e^{-\eta\,(\log N)^2},
    \end{align*}
    we obtain
    \begin{align*}
        \lim_{N\to\infty}\pk{\sum_{i=1}^m q_i(\mu_{N}(t_i)-\E{\mu_{N}(t_i)})<x}=\pk{\xi<x},
    \end{align*}
    implying that the vector $(\mu_{N}(t_1)-\E{\mu_{N}(t_1)},\ldots,\mu_{N}(t_m)-\E{\mu_{N}(t_m)})$ converges in distribution to a Gaussian vector, and moreover
    \begin{align}
        &\lim_{N\to\infty} \biggl(\mu_{N}(t_1)-\E{\mu_{N}(t_1)},\ldots,\mu_{N}(t_m)-\E{\mu_{N}(t_m)}\biggr)\notag\\
        &\qquad=\lim_{N\to\infty}\left(\frac{1}{N}\sum_{i,j=1}^{N}(a_{i,j}(t_1)-p(t_1)),\ldots\frac{1}{N}\sum_{i,j=1}^{N}(a_{i,j}(t_m)-p(t_m))\right). \label{a_vector}
    \end{align}
    Hence, it suffices to calculate the limiting covariance of the vector \eqref{a_vector}. Without any loss of generality we consider $m=2$ with  $t_1\leqslant t_2$. Then, as the processes $a_{i,j}(\cdot)$ are independent for different pairs $(i,j)$ with $i<j$,
    \begin{align*}
        &\operatorname{Cov}(\mu_{N}(t_1)-\E{\mu_{N}(t_1)},\mu_{N}(t_2)-\E{\mu_{N}(t_2)})\\
        &\qquad=\operatorname{Cov}\left(\frac{1}{N}\sum_{i,j=1}^{N}a_{i,j}(t_1),\frac{1}{N}\sum_{i,j=1}^{N}a_{i,j}(t_2)\right)\\
        &\qquad = \frac{4}{N^2}\sum_{\substack{i,j=1\\i<j}}^{N}\operatorname{Cov}\left(a_{i,j}(t_1),a_{i,j}(t_2)\right) + \frac{1}{N^2}\sum_{i=1}^{N}\operatorname{Cov}(a_{i,i}(t_1),a_{i,i}(t_2)).
    \end{align*}
    Recall from \eqref{eq:COV} that $\operatorname{Cov}(a_{i,j}(t_1),a_{i,j}(t_2))
        =p(t_1)(1-p(t_1))\,e^{-(\LON+\LOFF )(t_2-t_1)}$ for any $i,j\in[N]$ and $t_1\leqslant t_2$.
    Hence, the convergence of finite-dimensional distributions is proven, and the only remaining step to achieve the claim of \netheo{main} is to show tightness of the series $\mu_N(t) - \E{\mu_N(t)}$.
    
    To establish tightness, according to \cite[Theorem 13.5]{billingsley2013convergence}, or actually its `moment version' that applies under condition \cite[Eqn.\ (13.14)]{billingsley2013convergence}, it is sufficient to verify that for any $r\leqslant t$ and any $s\in[r,t]$ 
    \begin{align}
        \E{\abs{\frac{1}{N}\sum_{\substack{i,j=1\\ i\leqslant j}}^N X_{i,j}(r,s)}^{2}\abs{\frac{1}{N}\sum_{\substack{i,j=1\\ i\leqslant j}}^N X_{i,j}(s,t)}^{2}} \leqslant (F(t)-F(r))^2\label{tightness_claim}
    \end{align}
    for $F(t):=35\kappa   \,t$  with $\kappa = \max(\LON,\LOFF)$ and 
    \begin{align}
    X_{i,j}(s,t) := 2^{\mathbb{I}\{i\not= j\}}\bigl(a_{i,j}(t) - p(t) - a_{i,j}(s) + p(s)\bigr).\label{X_def}
    \end{align}
    So as to show \eqref{tightness_claim}, we expand the product on the left-hand side. Such a representation is useful because, using that the processes $X_{i,j}(\cdot,\cdot)$ are centered and independent, we can disregard any term of this expansion which contains some process $X_{i,j}(\cdot,\cdot)$ only once, as the contribution of such a term is zero. So, the only remaining terms contain either some process $X_{i,j}(\cdot,\cdot)$ four times for a given pair $(i,j)$, i.e., terms of the type
    \begin{align}
        X^2_{i,j}(r,s)\,X^2_{i,j}(s,t),\label{term_1}
    \end{align}
    or two process $X_{i,j}(\cdot,\cdot)$, $X_{k,\ell}(\cdot,\cdot)$ for different pairs $(i,j)$ and $(k,\ell)$, i.e., terms of the types
    \begin{align}
        &X^2_{i,j}(r,s)\,X_{k,\ell}^2(s,t),\label{term_2}\end{align}
        and
        \begin{align}
        &X_{i,j}(r,s)\,X_{k,\ell}(r,s)\,X_{i,j}(s,t)\,X_{k,\ell}(s,t).\label{term_3}
    \end{align}
    We proceed by separately assessing the contributions of terms of the types \eqref{term_1}, \eqref{term_2}, and \eqref{term_3}.

    \begin{itemize}
       \item[$\circ$] 
    The term \eqref{term_1} appears in the expansion of \eqref{tightness_claim} just once for each pair $(i,j)$ such that $i\leqslant j$. As a consequence, the factor in front of \eqref{term_1} is $1/N^4$. 
    \item[$\circ$] The term \eqref{term_2} also appears only once for each $i\leqslant j$, hence the factor in front of \eqref{term_2} is $1/N^4$ as well. 
    \item[$\circ$]
    The last term \eqref{term_3} is the most involved one. It appears four times for each $i\leqslant j$ and $k\leqslant \ell$, as $X_{i,j}(r,s)\,X_{k,\ell}(r,s)$ appears twice in the expansion of $(\sum_{i\leqslant j}X_{i,j}(r,s))^2$, and the same holds for $X_{i,j}(s,t)\,X_{k,\ell}(s,t)$. Hence, the factor in front of \eqref{term_3} is $4/N^4$. 
    \end{itemize}
    The above arguments entail that, recalling that $X_{i,i}(\cdot,\cdot)$ and $X_{i,j}(\cdot,\cdot)$ (with $i\not=j$) are treated differently in \eqref{X_def}, the left-hand side of \eqref{tightness_claim} can be rewritten as
    \begin{align}
         \frac{1}{N^4}\sum_{\substack{i,j=1 \\ i\leqslant j}}^{N}\E{X^2_{i,j}(r,s)\,X^2_{i,j}(s,t)}&+ \frac{1}{N^4}\sum_{\substack{i,j,k,\ell=1 \\ i\leqslant j,~k\leqslant \ell \\ (i,j)\not= (k,\ell)}}^{N}\E{X^2_{i,j}(r,s)\,X_{k,\ell}^2(s,t)}\notag\\
        &+ \frac{4}{N^4}\sum_{\substack{i,j,k,\ell=1 \\ i\leqslant j,~k\leqslant \ell \\ (i,j)\not= (k,\ell)}}^{N}\E{X_{i,j}(r,s)\,X_{k,\ell}(r,s)\,X_{i,j}(s,t)\,X_{k,\ell}(s,t)}.\notag
    \end{align}
    We continue by majorizing these three terms, as follows. For notational convenience, let the generic process $X(\cdot,\cdot)$ be distributed as $X_{1,1}(\cdot,\cdot)$, and, likewise, let $a(\cdot)$ be distributed as $a_{1,1}(\cdot)$. Recall that, by \eqref{X_def}, $X_{i,j}(\cdot,\cdot)$ has the same distribution as $2X(\cdot,\cdot)$ if $i\not=j$.
    It follows that the three terms are bounded from above by
    \begin{align}    \frac{16N(N+1)}{2N^4}\E{X^2(r,s)X^2(s,t)},\:\:&\:\:
         \frac{16N^2(N+1)^2}{4N^4}\E{X^2(r,s)}\E{X^2(s,t)} \notag
         \end{align}
         and
         \begin{align}\frac{64N^2(N+1)^2}{4N^4}\left(\E{X(r,s)\,X(s,t)}\right)^2, \notag
         \end{align}
         respectively.
         We thus find that for sufficiently large $N$  such that
         \begin{align*}
            \frac{N+1}{N^3}\leqslant 1,\qquad \frac{(N+1)^2}{N^2}\leqslant \frac{17}{16}\leqslant \frac{5}{4}, \qquad 
         \end{align*}the left-hand side of \eqref{tightness_claim} is majorized by 
         \begin{align}
        & 8\,\E{X^2(r,s)X^2(s,t)}+5\,\operatorname{Var}\left(X(r,s)\right)\operatorname{Var}\left(X(s,t)\right)+17\,\operatorname{Var}\left(X(r,s)\right)\operatorname{Var}\left(X(s,t)\right)\notag\\
        &\leqslant 8\,\E{X^2(r,s)X^2(s,t)}+22\,\operatorname{Var}\left(X(r,s)\right)\operatorname{Var}\left(X(s,t)\right). \label{eq:UB}
    \end{align}

The next step is to find upper bounds on each term given in \eqref{eq:UB}, i.e., $\E{X^2(r,s)X^2(s,t)}$ and $\operatorname{Var}\left(X(r,s)\right)\,\operatorname{Var}\left(X(s,t)\right)$. 
    Using that from \eqref{X_def} we have 
    \begin{align*}
        X(s,t)&=\begin{cases}
                p(s)-p(t)-1\quad &\text{with probability }\pk{a(t)=0,\,a(s)=1},\\
                p(s)-p(t) \quad &\text{with probability }\pk{a(t)=a(s)},\\
                p(s)-p(t)+1\quad &\text{with probability }\pk{a(t)=1,\,a(s)=0},
            \end{cases}
    \end{align*}
    we obtain, according to \eqref{P_def} and \eqref{p_def},
    \begin{align}
        \operatorname{Var}(X(s,t))&\leqslant \abs{p(s)-p(t)}^2 + 4\,\pk{a(t)\not = a(s)}\notag\\
        &= \abs{p(0)-\varrho^2}\left(e^{-(\LON+\LOFF ) t} - e^{-(\LON+\LOFF ) s}\right)^2\:+\notag\\
        &\quad\quad\quad 4\,\pk{\mathcal{P}_{a_{1,1}(0),1,1}(s) - \mathcal{P}_{a_{1,1}(0),1,1}(t)\text{ is odd}}\notag\\
        &\leqslant (1-e^{-2\,\kappa({t-s})})+ 4\,\pk{\mathcal{P}_{a_{1,1}(s),1,1}({t-s})\geqslant 1}\notag\\
        &=(1-e^{-2\kappa ({t-s})}) + 4\,(1-e^{-\kappa({t-s})})\leqslant 6\,\kappa ({t-s})\leqslant 6\,\kappa({t-r}).\label{tightness_2}
    \end{align}
    By a similar argumentation,
    \begin{align}
       \E{X^2(r,s)X^2(s,t)}&\leqslant 16\,\pk{\mathcal{P}_{a_{1,1}(s),1,1}({t-s})\geqslant 1}\pk{\mathcal{P}_{a_{1,1}(r),1,1}({r-s})\geqslant 1}\notag\\
       &\qquad + 4\abs{p(s)-p(t)}^2\pk{\mathcal{P}_{a_{1,1}(r),1,1}({r-s})\geqslant 1}\notag\\
       &\qquad  + 4\abs{p(r)-p(s)}^2\pk{\mathcal{P}_{a_{1,1}(s),1,1}({t-s})\geqslant 1}\notag\\
        &\qquad + \abs{p(s)-p(t)}^2\abs{p(r)-p(s)}^2\notag\\
        &\leqslant (16+8+8+16)\,\kappa^2\abs{t-s}\abs{r-s}\leqslant 48\,\kappa^2\abs{t-r}^2.\label{tightness_3}
    \end{align}
    Upon inserting the upper bounds \eqref{tightness_2} and \eqref{tightness_3} into \eqref{eq:UB}, we conclude
    \begin{align*}
        &\E{\abs{\frac{1}{N}\sum_{\substack{i,j=1\\ i\leqslant j}}^N X_{i,j}(r,s)}^{2}\abs{\frac{1}{N}\sum_{\substack{i,j=1\\ i\leqslant j}}^N X_{i,j}(s,t)}^{2}}\leqslant (384+792)\,\kappa^2\,\abs{t-r}^2\leq(F(t)-F(r))^2,
    \end{align*}
    recalling the definition of $F(\cdot)$. 
    This proves the tightness claim, by appealing to the sufficient condition \eqref{tightness_claim}.
    \QED

   \section{Proof of Theorem~\ref{true_eigenvalue_representation}}

    The main objective of this section is to prove Theorem~\ref{true_eigenvalue_representation}. 
    Following the ideas developed in \cite{10.1214/11-AOP734}, we normalize the adjacency matrix $A_N(t)$ in such a way that its entries have a constant variance $1/N$, and then isolate a centered stochastic part from this matrix. 
    In more concrete terms, we define the two auxiliary matrix-valued functions
    \begin{align}
        A^{\star}_{N}(t) &:= \frac{1}{\sqrt{Np(t)(1-p(t))}} A_N(t),\label{A*_def}\\
        H_N(t) & = \{h_{i,j}(t,N)\}_{i,j=1}^{N} := \left\{\frac{a_{i,j}(t)-p(t)}{\sqrt{Np(t)(1-p(t))}}\right\}_{i,j=1}^{N}.\label{H_def}
    \end{align}
Define the set of jump points for $H_N$ via
        \begin{align}
            \mathcal{T}(N):=\{t\in[0,T]\colon \exists\, i,j\in\{1,\ldots,N\}\colon h_{i,j}(t+,N)\not= h_{i,j}(t-,N) \} .\label{T_def}
        \end{align}
        The following proposition lists some important properties of the entries of the matrix $H_N(t)$. In particular, it provides an estimate on the jump points, which are at most distance $x$ apart. 
    \begin{prop}\label{model_basic_properties}
    $H_N(t)$ is random matrix with independent entries below the diagonal, satisfying
    \begin{align}
            \abs{h_{i,j}(t,N)}\leqslant \frac{C_1}{\sqrt{N}},\:\: \E{h_{i,j}(t,N)}=0,\:\:\E{\abs{h_{i,j}(t,N)}^2}=\frac{1}{N},\:\: \E{\abs{h_{i,j}(t,N)}^r}\leqslant\left(\frac{C_2}{\sqrt{N}}\right)^r\label{S2}
    \end{align}
    for some positive constants $C_1,C_2$. In particular, for all $t\in[0,T]$ we have almost surely
    \begin{align}
        \left\Vert H_N(t)\right\Vert \leqslant C_1\sqrt{N}.\label{H_full_bound}
    \end{align}
        In addition,
        for any $x\in[0,t]$ and sufficiently large $N$,
        \begin{align}
            \pk{\exists\,  t_1,t_2\in \mathcal{T}(N)\colon \abs{t_1-t_2}<x}\leqslant C N^4 x\label{A0}
        \end{align}
        for some positive constant $C$.
    \end{prop}
    In fact, \eqref{S2} can be regarded as the uniform analogue of the moment bounds \cite[Eqn.\ (3.2)]{10.1214/11-AOP734}, and \eqref{A0} will be instrumental in finding an event $\Omega_N$ in \netheo{true_eigenvalue_representation} independent of $t$. 

\medskip

    \proofprop{model_basic_properties} 
Using \eqref{H_def} we obtain that
\begin{align*}
\abs{h_{i,j}(t,N)}=\frac{\abs{a_{i,j}(t)-p(t)}}{\sqrt{Np(t)(1-p(t))}}\leqslant \frac{1}{\sqrt{p^-(1-p^+)}}\frac{1}{\sqrt{N}}.
\end{align*}
To obtain \eqref{S2}, observe that
\begin{align*}
    \E{h_{i,j}(t,N)}&=\frac{\E{a_{i,j}(t)}-p(t)}{\sqrt{Np(t)(1-p(t))}}=0,\\
    \E{\abs{h_{i,j}(t,N)}^2}&=\frac{\operatorname{Var}\bigl(a_{i,j}(t)\bigr)}{Np(t)(1-p(t))}=\frac{p(t)(1-p(t))}{Np(t)(1-p(t))}=\frac{1}{N},\\
    \E{\abs{h_{i,j}(t,N)}^r}&=\left(\frac{\E{\abs{a_{i,j}(t)-p(t)}}}{\sqrt{Np(t)(1-p(t))}}\right)^r\leqslant \left(\frac{1}{\sqrt{p^-(1-p^+)}}\frac{1}{\sqrt{N}}\right)^r.
\end{align*}

To show \eqref{H_full_bound} we use the first statement of \eqref{S2} obtaining that
\begin{align*}
   \left\Vert H_N(t)\right\Vert \leqslant \max_{i\in[N]}\sum_{j\in [N]}\abs{h_{i,j}(t,G)} \leqslant  \max_{i\in[N]}\sum_{j\in [N]}\frac{C_1}{\sqrt{N}} = C_1\sqrt{N}.
\end{align*}

Finally, to show \eqref{A0}, let us consider only $x<T$ and define for any $a\in\{0,1\}$ and vector ${\boldsymbol t}\equiv (t_1,t_2)\in[0,T]^2$, 
\[{\mathcal E}_{i,j}\big[{a},{\boldsymbol t}\big]
:= \abs{\mathcal{P}_{{a},i,j}(t_1)-\mathcal{P}_{{a},i,j}(t_2)},
\]
where recall $\mathcal{P}_{{a}, i, j}$ is defined in \eqref{P_def}. Also define a set ${\mathcal T}_x:=\{{\boldsymbol t}\in[0,T]^2: \abs{t_1-t_2}<x\}$. The event ${\mathcal E}_{i,j}({a},{\boldsymbol t})$ tells us how many times the edge between vertices $i$ and $j$ changed state between times $t_1$ and $t_2$, given that at time $t=0$ it was in state $a$.
Thus, the probability in the right-hand side of \eqref{A0} can be bounded as follows:
\begin{align*}
    &\pk{\exists\,  {\boldsymbol t}\in\mathcal{T}(N) \colon \abs{t_1-t_2}<x}
    \leqslant \pk{\exists\,  {\boldsymbol t}\in {\mathcal T}_x,\, i,j\in\{1,\ldots,N\} \colon {\mathcal E}_{i,j}\big[a_{i,j}(0),{\boldsymbol t}\big]\geqslant 2}\\
&\qquad+ \pk{\exists\,  {\boldsymbol t}\in {\mathcal T}_x, {\boldsymbol i},{\boldsymbol j}\in\{1,\ldots,N\}^2 \colon {\boldsymbol i}\not={\boldsymbol j}, ~
\,{\mathcal E}_{i_1,j_1}\big[a_{i_1,j_1}(0),{\boldsymbol t}\big]\geqslant 1, 
\,{\mathcal E}_{i_2,j_2}\big[a_{i_2,j_2}(0),{\boldsymbol t}\big]\geqslant 1
}\\
    &\quad\leqslant N^2\,\pk{\exists\, {\boldsymbol t}\in {\mathcal T}_x \colon {\mathcal E}_{1,1}\big[a_{1,1}(0),{\boldsymbol t}\big]\geqslant 2}\\
    &\qquad+ N^4\,\pk{\exists\,  {\boldsymbol t}\in{\mathcal T}_x:~\,{\mathcal E}_{1,1}\big[a_{1,1}(0),{\boldsymbol t}\big]\geqslant 1, 
\,{\mathcal E}_{2,2}\big[a_{2,2}(0),{\boldsymbol t}\big]\geqslant 1},
\end{align*}
with ${\boldsymbol i}\equiv(i_1,i_2)$ and ${\boldsymbol j}\equiv(j_1,j_2)$.
From \eqref{P_def} we obtain further that, with $\xi^1$ distributed as an exponential random variable with parameter $\LON$ and $\xi^0$ distributed as an exponential random variable with parameter $\LOFF$,
\begin{align}
   &\pk{\exists \,  {\boldsymbol t}\in  {\mathcal T}_x   \colon {\mathcal E}_{1,1}\big[a_{1,1}(0),{\boldsymbol t}\big]\geqslant 2} \notag\\
   &\qquad \leqslant \pk{\exists\,  k\in\{0,\ldots,\lceil T/x\rceil\}\colon {\mathcal E}_{1,1}\big[{a_{1,1}(0)},{(kx+3x,kx)}\big] \geqslant 2}\notag\\
   &\qquad  \leqslant \sum_{k=0}^{\lceil T/x\rceil}\pk{{\mathcal E}_{1,1}\big[{a_{1,1}(0)},{(kx+3x,kx)}\big] \geqslant 2},\label{TERM1}
\end{align}
by the union bound. With $w_a(x):= \pk{{\mathcal E}_{1,1}\big[{a},{(x,0)}\big] \geqslant 2}$ for $a\in\{0,1\}$, using the Markov property,
\[\pk{{\mathcal E}_{1,1}\big[{a_{1,1}(0)},{(kx+3x,kx)}\big] \geqslant 2}\leqslant\max\{w_0(3x),w_1(3x)\}.\]
It is a matter of straightforward calculus to verify that 
\begin{align*}
w_1(x)&= 1 - e^{-\LON  x} - \LON  \left(\frac{e^{-\LOFF  x}-e^{-\LON  x}}{\LON -\LOFF }\right)=x^2 \frac{\LON \LOFF }{2}+O(x^3),
\end{align*}
and the same expression for $w_0(x)$ but then with the roles of $\LON $ and $\LOFF $ being swapped (where it is easily verified that this bound also holds in case $\LON =\LOFF $ by performing a limiting argument).

Hence,  there is an $\varepsilon>0$ so that for all $x\in[0,t]$, \eqref{TERM1} is majorized by
\[\left(\frac{T}{x}+2\right)\left(\frac{\LON\LOFF }{2}+\varepsilon\right)\cdot 9\,x^2. \]
Similarly, there exists an $\varepsilon'>0$ such that for all $x\in[0,t]$
\begin{align*}
    &\pk{\exists \, {\boldsymbol t}\in{\mathcal T}_x:~\,{\mathcal E}_{1,1}\big[a_{1,1}(0),{\boldsymbol t}\big]\geqslant 1, 
\,{\mathcal E}_{2,2}\big[a_{2,2}(0),{\boldsymbol t}\big]\geqslant 1}\\
    &\qquad\leqslant \pk{\exists\,  k\in\{0,\ldots,\lfloor T/x\rfloor\}\colon {\mathcal E}_{1,1}\big[a_{1,1}(0),(kx+3x,k x)\big]\geqslant 1, 
\,{\mathcal E}_{2,2}\big[a_{2,2}(0),(kx+3x,kx)\big]\geqslant 1}\\
    &\qquad\leqslant \left( \frac{T}{x}+2\right) \left(\max\big\{\pk{\xi^0\leqslant 3x},\pk{\xi^1\leqslant 3x}\big\}\right)^2\leqslant \left( \frac{T}{x}+2\right)(\kappa+\varepsilon')^2\cdot9\,x^2.
\end{align*}
Hence, the claim \eqref{A0} follows.
\QED

\medskip

    To make our notation more compact, we introduce the auxiliary function
    \begin{align*}
        q(t) := \frac{p(t)}{1-p(t)}.
    \end{align*}
    Using \eqref{p_bounds}, for all $t\in[0,T]$,
    \begin{align*}
        0<q^-:=\frac{p^-}{1-p^-}\leqslant q(t)\leqslant  \frac{p^+}{1-p^+}=:q^+ < \infty.
    \end{align*}
    
    The matrices $A^{\star}_N(t)$ and $H_N(t)$ can be expressed in one another via
    \begin{align}
        A^{\star}_N(t) = H_N(t) + \sqrt{Nq(t)}\,\mathcal{E}_N,\label{A_in_H}
    \end{align}
    where $\mathcal{E}_N$ is an $N\times N$ matrix of which all entries equal $1/N$.
    Define by $\mu_{N}^{\star}(t)$ the principal eigenvalue of the matrix $A^{\star}_{N}(t)$. Then, according to \eqref{A*_def},
    \begin{align}
    \mu_{N}(t) = \sqrt{Np(t)(1-p(t))}\,\mu^{\star}_{N}(t),\label{mu_from_mu*}
    \end{align}
    so that, combining \eqref{mu_from_mu*} with \eqref{H_def}, \netheo{true_eigenvalue_representation} is an immediate corollary of the following lemma.

    \begin{lem}\label{mu_sum_expression}
    There exist a set of events $\Omega_N$ such that
    \begin{align*}
        \pk{\Omega_N^{\rm c}}\leqslant e^{-\nu(\log N)^2}
    \end{align*}
    and on $\Omega_N$
    \begin{align*}
        \mu^{\star}_{N}(t) = \E{\mu^{\star}_{N}(t)} + \frac{1}{N}\sum_{i,j\in\{1,\ldots,N\}}h_{i,j}(t,N) + \frac{\mathcal{G}_N(t)}{\sqrt{p(t)(1-p(t))}}\left(\frac{(\log N)^{4}}{N}\right),
    \end{align*}
    where $\limsup_{N\to\infty}\sup_{t\in[0,T]} \abs{\mathcal{G}_{N}(t)}<\infty$. 
\end{lem}

To show \nelem{mu_sum_expression} we use the following results.  From this point, let $\vk e_N\in \R^N$ be a unit vector, i.e., each entry of this vector is equal to $1/\sqrt{N}$. As the dimension of $\vk e_N$ is usually clear from the context, we omit the subscript and write simply $\vk e$. The main tool when analyzing the eigenvalue $\mu^{\star}_{N}$ is the following series expansion.

\begin{lem}\label{lem_mu_recurrent_expansion}
    On the set of events $\Omega_N$ defined in \nelem{uniform_event},  the principal eigenvalue $\mu^{\star}_{N}(t)$ of the matrix $A^{\star}_{N}(t)$ satisfies the following equation:
    \begin{align}
        \mu^{\star}_{N}(t)=\sqrt{Nq(t)}\sum_{k= 0}^\infty\sprod{\vk e, \left(\frac{H_{N}(t)}{\mu^{\star}_{N}(t)}\right)^k\vk e}.\label{mu_recurrent_expansion}
    \end{align}
\end{lem}

In our proof we use that $\mu_N^{\star}(t)$ is asymptotically equivalent to a partial sum of the series in the right hand side of \eqref{mu_recurrent_expansion}. To verify this, we need to show that the tail of the series is, as $N\to\infty$, asymptotically negligible with respect to that partial sum. Hence, taking into account the inequality
\begin{align*}
    \left\vert\sprod{\vk e, H^k_N(t)\,\vk e}\right\vert\leqslant \Vert H^k_{N}(t)\Vert \leqslant \left\Vert H_{N}(t)\right\Vert^k,
\end{align*}
we should derive an upper bound (uniform in $t\in[0,T]$) on $\left\Vert H_{N}(t)\right\Vert$, as well as a lower bound (uniform in $t\in[0,T]$) on $\mu_N^{\star}(t)$, such that they allow us to verify that, as $N\to\infty$,
\begin{align}
    \frac{\left\Vert H_{N}(t)\right\Vert}{\mu_N^{\star}(t)}\to 0\label{tail_claim}
\end{align}
uniformly for $t\in[0,T]$. Such bounds, however, do not exist in the almost sure sense, but we remedy this by finding bounds that apply on events $\Omega_N$ of sufficiently high probability. Consider first the bound for $\left\Vert H_{N}(t)\right\Vert$. The following lemma plays a key role in it.

\begin{lem}\label{uniform_event} There exists a sequence of events $\Omega_N$ such that for sufficiently large $N$ 
\begin{align}
    \pk{\Omega_N^{\rm c}}\leqslant e^{-\nu(\log N)^2},\label{Omega_bound}
\end{align}
    and the two following bounds hold on $\Omega_N$, for sufficiently large $N$, for all $t\in[0,T]$ and all $k\in\{1,\ldots, \lceil\log N\rceil\}$:
    \begin{align}
        &\left\Vert H_N(t)\right\Vert\leqslant 2 + \frac{(\log N)^{2}}{N^{1/4}},\label{A1}\\
        &\abs{\sprod{\vk e,H^k_N(t)\,\vk e}-\E{\sprod{\vk e,H^k_N(t)\,\vk e}}}<\mathfrak{C}^{k}\frac{(\log N)^{2k}}{\sqrt{N}}\label{A2}
    \end{align}
    for some positive constant $\mathfrak{C}$.
\end{lem}

For each particular value of $t\in[0,T]$ the bounds \eqref{A1} and \eqref{A2} have already been established; see Lemma 4.3 and Lemma 6.5 (combined with Eqn.~(6.20)) in  \cite{10.1214/11-AOP734}. However, in the proof of the {\it functional} convergence it is crucial in \nelem{uniform_event} that the event $\Omega_N$ can be chosen {\it uniformly} for all $t\in[0,T]$, and this property cannot be straightforwardly extracted from the results presented in \cite{10.1214/11-AOP734}. This entails that a crucial element in the proof of Lemma \ref{uniform_event} amounts to appropriately choosing the set $\Omega_N$.

With \eqref{A1} controlling $\Vert H_N(t)\Vert$, the following lemma provides us with the desired uniform bound on $\mu_N^{\star}(t)$, so that we can establish \eqref{tail_claim}.

\begin{lem}\label{mu_full_expression}
On $\Omega_N$ defined in \nelem{uniform_event}, for all $t\in[0,T]$,
\begin{align}
    \mu^{\star}_{N}(t)=\sqrt{Nq(t)} + \tilde{\mathcal{G}}_N(t)\frac{(\log N)^{2}}{\sqrt{N}},\label{mu_precise_expression}
\end{align}
for some $\tilde{\mathcal{G}}_{N}(t)$ such that $\limsup_{N\to\infty}\sup_{t\in[0,T]}\abs{\tilde{\mathcal{G}}_{N}(t)}$ is finite almost surely.
\end{lem}

In fact, it is directly seen that this approximation is a truncation of \eqref{mu_recurrent_expansion} at the first term. Again, compare it with the counterpart for a given value of $t$, established in \cite[Eqn.\ (6.8)]{10.1214/11-AOP734}; the uniform bound on the residual term $\tilde{\mathcal{G}}_{N}(t)$ is needed in the proof of the functional convergence and cannot be straightforwardly obtained from \cite{10.1214/11-AOP734}.
The proofs of Lemmas  \ref{lem_mu_recurrent_expansion}, \ref{uniform_event}, and \ref{mu_full_expression} can be found in the appendix. In the remainder of this section, we prove Lemma \ref{mu_sum_expression}, using Lemmas \ref{lem_mu_recurrent_expansion}, \ref{uniform_event}, and \ref{mu_full_expression}.

\medskip

\prooflem{mu_sum_expression} We mimic the line of reasoning of the proof of \cite[Eqn (6.8)]{10.1214/11-AOP734}. As we need to obtain a uniformly bounded function $\mathcal{G}_N(t)$, however, we need to provide bounds on various terms in the series expansion \eqref{mu_recurrent_expansion}.

As we are going to use the bounds established in \nelem{uniform_event}, we shall use the event $\Omega_N$ defined there.
In this proof we are successively going through the following two steps. In the first place, we truncate the series that appears in the statement of \nelem{lem_mu_recurrent_expansion} and replace the inner products $\sprod{\vk e, H^k_N(t)\,\vk e}$ in the numerator appearing in \eqref{mu_recurrent_expansion} by its expectation. Second, we analyze the difference between the solution of the resulting equation (which is not stochastic in nature) and the solution of the initial equation (which is $\mu^{\star}_{N}(t)$). 

\medskip 

\noindent {\sc Step 1.}
We set the truncation level at $[\log N]$. To bound the tail of the series we introduce an upper bound on each of its terms. Appealing to \nelem{mu_full_expression}, we obtain that, on $\Omega_N$ and for sufficiently large $N$,
\begin{align}
    \frac{\sqrt{Nq^-}}{2}\leqslant\mu^{\star}_N(t)\leqslant 2\sqrt{Nq^+}\label{mu_rough_bound}.
\end{align}
Combining  this with the upper bound on the norm of $H_N(t)$ stated in \eqref{A1} (in \nelem{uniform_event}), we obtain that, for $N$ large enough, uniformly over $\Omega_N$,
\begin{align}
    \frac{\left\Vert H_N(t) \right\Vert}{\mu^{\star}_{N}(t)} \leqslant \frac{6}{\sqrt{Nq^-}}\label{H/mu_bound}.
\end{align}
We rewrite the infinite series \eqref{mu_recurrent_expansion} as follows:
\begin{align}
    \mu^{\star}_{N}(t) &=\sqrt{Nq(t)}\sum_{k=0}^{[\log N]}\frac{\sprod{\vk e,H^k_{N}(t)\,\vk e}}{\bigl(\mu^{\star}_{N}(t)\bigr)^k} + \sqrt{Nq(t)}\sum_{k=[\log N]+1}^{\infty}\frac{\sprod{\vk e,H^k_N(t)\,\vk e}}{\bigl(\mu^{\star}_{N}(t)\bigr)^k}\notag\\
    &=\sqrt{Nq(t)}\sum_{k=0}^{[\log N]}\frac{\sprod{\vk e,H^k_{N}(t)\,\vk e}}{\bigl(\mu^{\star}_{N}(t)\bigr)^k} + E_{1,N}(t)\,N^{\frac{2-[\log N]}{4}}\notag\\
    &=\sqrt{Nq(t)}\sum_{k=0}^{[\log N]}\frac{\E{\sprod{\vk e,H^k_{N}(t)\,\vk e}}}{\bigl(\mu^{\star}_{N}(t)\bigr)^k} +\frac{\sqrt{Nq(t)}}{\mu^{\star}_{N}(t)}\sprod{\vk e,H_N(t)\,\vk e}\notag\\
    &\qquad+E_{2,N}(t)\frac{(\log N)^{4}}{N}+ E_{1,N}(t)\,N^{\frac{2-[\log N]}{4}}\notag\\
    &=\sqrt{Nq(t)}\sum_{k=0}^{[\log N]}\frac{\E{\sprod{\vk e,H^k_{N}(t)\,\vk e}}}{\bigl(\mu^{\star}_{N}(t)\bigr)^k} +\frac{\sqrt{Nq(t)}}{\mu^{\star}_{N}(t)}\sprod{\vk e,H_N(t)\,\vk e}\label{mu*_in_E}+ E_{3,N}(t)\frac{(\log N)^{4}}{N},
\end{align}
using the notation
\begin{align*}
E_{1,N}(t)&:=N^{\frac{[\log N]-2}{4}}\sqrt{Nq(t)}\sum_{k=[\log N]+1}^{\infty}\frac{\sprod{\vk e,H^k_N(t)\,\vk e}}{\bigl(\mu^{\star}_{N}(t)\bigr)^k},\\
E_{2,N}(t)&:=\frac{N}{(\log N)^4}\sqrt{Nq(t)}\sum_{k=2}^{[\log N]}\left(\frac{\sprod{\vk e,H^k_N(t)\,\vk e} - \E{\sprod{\vk e,H^k_{N}(t)\,\vk e}}}{\bigl(\mu^{\star}_{N}(t)\bigr)^k}\right)\\
E_{3,N}(t)&:=\frac{N^{\frac{6-[\log N]}{4}}}{(\log N)^4}E_{1,N}(t) + E_{2,N}(t).
\end{align*}
Then we bound $E_{1,N}(t)$, $E_{2,N}(t)$, and $E_{3,N}(t)$ by constants, relying on inequalities derived earlier in this proof. By \eqref{H/mu_bound}, on $\Omega_N$ and for sufficiently large $N$, we can majorize $\abs{E_{1,N}(t)}$ by 
\begin{align}
    \abs{E_{1,N}(t)}&\leqslant \frac{\sqrt{N}}{N^{\frac{2-[\log N]}{4}}}\sqrt{q^+} \sum_{k=[\log N]}^{\infty}\left(\frac{6}{\sqrt{Nq^-}}\right)^k\notag\\
    &\leqslant \frac{\sqrt{N}}{N^{\frac{2-[\log N]}{4}}}\sqrt{q^+} \left(\frac{6}{\sqrt{Nq^-}}\right)^{[\log N]} \sum_{k=0}^{\infty}\left(\frac{6}{\sqrt{Nq^-}}\right)^k\notag\\
    &=\sqrt{q^+}\left(\frac{36}{q^-}\right)^{\frac{[\log N]}{2}}N^{\frac{1-[\log N]}{2}}N^{\frac{[\log N]-2}{4}}\frac{1}{1-\frac{6}{\sqrt{Nq^-}}}\notag\\
   &\leqslant 2\sqrt{q^+}\left(\frac{36}{q^-}\right)^{\frac{[\log N]}{2}}N^{-\frac{[\log N]}{4}}
  = 2\sqrt{q^+}\left(\frac{1296}{N(q^-)^2}\right)^{\frac{[\log N]}{4}}
  \:\leqslant\: 2\sqrt{q^+}.
  \label{E1_bound}
\end{align}
By \eqref{mu_rough_bound}, in combination with \eqref{A2} (from \nelem{uniform_event}), on $\Omega_N$ and for sufficiently large $N$ ,
\begin{align}
    \abs{E_{2,N}(t)}&\leqslant \frac{N}{(\log N)^4}\sqrt{Nq^+}\sum_{k=2}^{[\log N]}\frac{2^k\mathfrak{C}^{k}(\log N)^{2k}}{\sqrt{N}\left(\sqrt{Nq^-}\right)^k}\notag\\
    &= \frac{4\mathfrak{C}^2\sqrt{q^+}}{q^-}\sum_{k=0}^{[\log N]-2}\frac{2^k\mathfrak{C}^{k}(\log N)^{2k}}{\left(\sqrt{Nq^-}\right)^k}\notag\\
    &\leqslant\frac{4\mathfrak{C}^2\sqrt{q^+}}{q^-}\sum_{k=0}^{[\log N]-2}\left(\frac{2\mathfrak{C}(\log N)^2}{\sqrt{Nq^-}}\right)^k\leqslant \frac{4\mathfrak{C}^2\sqrt{q^+}}{q^-}\sum_{k=0}^{[\log N]-2}\frac{1}{2^k}\leqslant \frac{8\mathfrak{C}^2\sqrt{q^+}}{q^-}.\label{E2_bound}
\end{align}
Hence, upon combining the upper bound \eqref{E1_bound} and \eqref{E2_bound}, we find that the remaining term in \eqref{mu*_in_E} can be majorized as follows:
\begin{align}
    \abs{E_{3,N}(t)}&\leqslant 2\sqrt{q^+} + \frac{8\mathfrak{C}^2\sqrt{q^+}}{q^-}.\label{E3_bound}
\end{align}
This completes the first step of the proof.

\medskip 

\noindent {\sc Step 2.}
Define $\bar{\mu}_N(t)$ as a solution for $\mu$ in the equation
\begin{align}
\mu=\sqrt{Nq(t)}\sum_{k=0}^{[\log N]}\frac{\E{\sprod{\vk e,H^k_{N}(t)\,\vk e}}}{\mu^k}\label{mu_bar_def}
\end{align}
for 
 $\mu \geqslant \sqrt{Nq(t)}.$
One readily observes that both sides of \eqref{mu_bar_def} are continuous in $\mu$, the left-hand side is increasing in $\mu$, whereas the right-hand side is decreasing in $\mu$, thus leading to a unique solution. Moreover, as $\E{\sprod{\vk e,H_{N}(t)\,\vk e}} = 0$ and $\E{\sprod{\vk e,H^2_{N}(t)\,\vk e}}>0$, for large enough 
$N\in\N$
\begin{align*}
\sqrt{Nq(t)}\leqslant\sqrt{Nq(t)}\sum_{k=0}^{[\log N]}\frac{\E{\sprod{\vk e,H^k_{N}(t)\,\vk e}}}{\bigl(Nq(t)\bigr)^{k/2}},
\end{align*}
and 
\begin{align*}
    \lim_{\mu\to\infty}\sqrt{Nq(t)}\sum_{k=0}^{[\log N]}\frac{\E{\sprod{\vk e,H^k_{N}(t)\,\vk e}}}{\mu^k} = \sqrt{Nq(t)}<\infty.
\end{align*}
Hence, the solution of \eqref{mu_bar_def} exists and is unique. Similarly to what is done in the proof of \nelem{mu_full_expression}, it can be shown that one can choose $\bar{\mu}_N(t)>7$ such that
\begin{align*}
    \bar{\mu}_N(t) = \sqrt{Nq(t)} + \tilde{g}_N(t),
\end{align*}
where $\tilde{g}_N(t)$ is such that $\lim_{N\to\infty}\sup_{t\in[0,T]}\abs{\tilde{g}_N(t)}=0$. In particular, it means that $\bar{\mu}_N(t)$ satisfies the inequality 
\begin{align}
    \frac{\sqrt{Nq^-}}{2}\leqslant\bar{\mu}_N(t)\leqslant 2\sqrt{Nq^+}\label{bar_mu_rough_bound}
\end{align}
(being similar to the inequality \eqref{mu_rough_bound} that we encountered above). 
It is crucial here that $\bar{\mu}_N(t)$ is deterministic, meaning that it appears in $\E{\mu^{\star}_{N}(t)}$ as well. So, for the assertion claimed in \nelem{mu_sum_expression}, the only important part we need to analyze is the difference
\begin{align}
\zeta_N(t):=\mu^{\star}_{N}(t) - \bar{\mu}_N(t).\label{zeta_def}
\end{align}
From \nelem{mu_full_expression} it follows that
\begin{align}
    \lim_{N\to\infty}\sup_{t\in[0,T]}\abs{\zeta_N(t)}=0\label{zeta_to_zero}
\end{align}
on $\Omega_N$, but this property is not sufficiently sharp in order to prove the claimed assertion. So as to analyze the function $\zeta_N(t)$ more precisely, we combine the series \eqref{mu*_in_E} and \eqref{mu_bar_def}, resulting in
\begin{align}
\zeta_N(t) &= \frac{\sqrt{Nq(t)}}{\mu^{\star}_{N}(t)}\sprod{\vk e,H_N(t)\,\vk e} \notag\\
&\qquad+ \sqrt{Nq(t)}\sum_{k=0}^{[\log N]}\E{\sprod{\vk e,H^k_{N}(t)\,\vk e}}\left(\frac{1}{\bigl(\mu^{\star}_{N}(t)\bigr)^k} - \frac{1}{\bar{\mu}^k_N(t)}\right)+ E_{3,N}(t)\frac{(\log N)^{4}}{N}\notag\\
&= \sprod{\vk e,H_N(t)\,\vk e} - \frac{\mu^{\star}_{N}(t) - \sqrt{Nq(t)}}{\mu^{\star}_{N}(t)}\sprod{\vk e, H_N(t)\,\vk e}\notag\\
&\qquad+ \sqrt{Nq(t)}\sum_{k=2}^{[\log N]}\frac{\E{\sprod{\vk e,H^k_{N}(t)\,\vk e}}}{\bigl(\mu_{N}^{\star}(t)\bigr)^k}\left(\frac{\bar{\mu}^k_N(t)-(\mu_{N}^{\star}(t))^k}{\bar{\mu}^k_N(t)}\right)+ E_{3,N}(t)\frac{(\log N)^{4}}{N}\notag.
\end{align}
This can be further rewritten as
\begin{align}
\zeta_N(t)&= \sprod{\vk e,H_N(t)\,\vk e} - \frac{\tilde{\mathcal{G}}_N(t)\frac{(\log N)^{2}}{\sqrt{N}}}{\mu^{\star}_{N}(t)}\sprod{\vk e, H_N(t)\,\vk e}+ E_{4,N}(t)\frac{(\log N)^4}{N}+ E_{3,N}(t)\frac{(\log N)^{4}}{N}\notag\\
&= \sprod{\vk e,H_N(t)\,\vk e} + \big(E_{5,N}(t) + E_{4,N}(t)+E_{3,N}(t)\big)\frac{(\log N)^{4}}{N}\notag\\
&= \sprod{\vk e,H_N(t)\,\vk e} + E_{6,N}(t)\frac{(\log N)^{4}}{N}\label{zeta_last},
\end{align}
where the functions $E_{4,N}(t)$, $E_{5,N}(t)$, and $E_{6,N}(t)$ are given by
\begin{align*}
    E_{4,N}(t)&:=\frac{N}{(\log N)^4}\sqrt{Nq(t)}\sum_{k=2}^{[\log N]}\frac{\E{\sprod{\vk e,H^k_{N}(t)\,\vk e}}}{\bigl(\mu_{N}^{\star}(t)\bigr)^k}\left(\frac{\bar{\mu}^k_N(t)-\bigl(\mu^{\star}_{N}(t)\bigr)^k}{\bar{\mu}^k_N(t)}\right),\\
    E_{5,N}(t)&:=-\frac{N}{(\log N)^4}\frac{\tilde{\mathcal{G}}_N(t)\frac{(\log N)^{2}}{\sqrt{N}}}{\mu^{\star}_{N}(t)}\sprod{\vk e, H_N(t)\,\vk e},\\
    E_{6,N}(t)&:=E_{3,N}(t) + E_{4,N}(t) + E_{5,N}(t).
\end{align*}
This means that the next goal is to bound $E_{6,N}(t)$. We treat the different parts of $E_{6,N}(t)$ separately. Combining \eqref{H/mu_bound} with \eqref{mu_rough_bound} and \eqref{H_full_bound} we can show that, on $\Omega_N$ with sufficiently large $N$ and $k\leqslant \log N$, 
\begin{align}
    \abs{\frac{\E{\sprod{\vk e,H^k_{N}(t)\,\vk e}}}{\bigl(\mu_{N}^{\star}(t)\bigr)^k}} &\leqslant \abs{\frac{\E{\sprod{\vk e,H^k_{N}(t)\,\vk e}\mathbb{I}_{\Omega_N}}}{\bigl(\mu_{N}^{\star}(t)\bigr)^k}} + \abs{\frac{\E{\sprod{\vk e,H^k_{N}(t)\,\vk e}\mathbb{I}_{\Omega_N^{\rm c}}}}{\bigl(\mu_{N}^{\star}(t)\bigr)^k}}  \notag\\
    &\leqslant \abs{\frac{\sup_{\Omega_N}\sprod{\vk e,H^k_{N}(t)\,\vk e}}{\bigl(\mu_{N}^{\star}(t)\bigr)^k}} + \abs{\frac{\left(C_1\sqrt{N}\right)^{k}\pk{\Omega_N^{\rm c}}}{\bigl(\sqrt{Nq^{-}}/2\bigr)^k}}\notag\\
    &\leqslant \left(\frac{6}{\sqrt{Nq^-}}\right)^k + \left(\frac{2C_1}{\sqrt{q^-}}\right)^ke^{-\nu(\log N)^2} \notag\\
    &\leqslant \frac{36}{Nq^{-}}\left(\frac{6}{\sqrt{Nq^{-}}}\right)^{k-2} + \frac{4C_1^2}{Nq^{-}}\left(\frac{2C_1}{\sqrt{q^{-}}} e^{-\frac{\nu(\log N)^2 - \log N}{k-2}}\right)^{k-2}\notag\\
    &\leqslant \frac{36}{Nq^{-}}\left(\frac{6}{\sqrt{Nq^-}}\right)^{k-2} + \frac{4C_1^2}{Nq^{-}}\left(\frac{2C_1}{\sqrt{q^{-}}} e^{-\frac{\frac{\nu}{2}(\log N)^2}{k}}\right)^{k-2}\notag\\
    &\leqslant \frac{36}{Nq^{-}}\left(\frac{6}{\sqrt{Nq^-}}\right)^{k-2} + \frac{4C_1^2}{Nq^{-}}\left(\frac{2C_1}{\sqrt{q^{-}}} e^{-\frac{\nu\log N}{2}}\right)^{k-2}\notag\\
    &\leqslant \frac{1}{N}\frac{4C_1^2+36}{q^-}\left(\frac{2C_1 + 6}{\sqrt{q^{-}}}N^{-\frac{\min\{1,\nu\}}{2}}\right)^{k-2}.\label{E4.1}
\end{align}

Combining \eqref{mu_rough_bound} and \eqref{bar_mu_rough_bound}, for any $m\in[k-1]$
\begin{align}
    \abs{\bar{\mu}_N^m(t)\bigl(\mu^{\star}_{N}(t)\bigr)^{k-1-m}}\leqslant\left(2\sqrt{Nq^+}\right)^{k-1}.\label{mu_prod_upper_bound}
\end{align}
Hence, using \eqref{mu_prod_upper_bound} together with \eqref{bar_mu_rough_bound} and \eqref{zeta_to_zero}, we obtain for sufficiently large $N$  that
\begin{align}
    \abs{\frac{\bar{\mu}^k_N(t)-\bigl(\mu^{\star}_{N}(t)\bigr)^k}{\bar{\mu}^k_N(t)}} &= \frac{\abs{\bar{\mu}_N(t)-\mu^{\star}_{N}(t)}\abs{\sum\limits_{m=0}^{k-1}\bar{\mu}_N^m(t)\bigl(\mu^{\star}_{N}(t)\bigr)^{k-1-m}}}{\abs{\bar{\mu}^k_N(t)}}\notag\\
    &\leqslant \frac{\abs{\zeta_N(t)}k\left(2\sqrt{Nq^+}\right)^{k-1}}{\left(\frac{1}{2}\sqrt{Nq^-}/2\right)^k}\leqslant \frac{2k\left(2\sqrt{Nq^+}\right)^{k-1}}{\left(\frac{1}{2}\sqrt{Nq^-}\right)^k}\label{E4.2}
\end{align}
By combining \eqref{E4.1} and \eqref{E4.2}, on $\Omega_N$ with sufficiently large $N$, the term $E_{4,N}(t)$ can be bounded as follows:
\begin{align}
\abs{E_{4,N}(t)}&\leqslant \frac{4C_1^2 + 36}{q^{-}(\log N)^4}\sqrt{Nq^+}\sum_{k=2}^{[\log N]}\left(\frac{2C_1 + 6}{\sqrt{q^{-}}N^{\frac{\min\{1,\nu\}}{2}}}\right)^{k-2}\left(\frac{2k\left(2\sqrt{Nq^+}\right)^{k-1}}{\left(\frac{1}{2}\sqrt{Nq^-}\right)^k}\right)\notag\\
&= \frac{2C_1^2 + 18}{q^{-}(\log N)^4}\sum_{k=2}^{[\log N]}\left(\frac{2C_1 + 6}{\sqrt{q^{-}}N^{\frac{\min\{1,\nu\}}{2}}}\right)^{k-2}\left(\frac{\left(\sqrt[k]{k}\right)^k\left(2\sqrt{Nq^+}\right)^{k}}{\left(\frac{1}{2}\sqrt{Nq^-}\right)^k}\right)\notag\\
&= \frac{2C_1^2 + 18}{q^-(\log N)^4}\sum_{k=2}^{[\log N]}\left(\frac{2C_1 + 6}{\sqrt{q^{-}}N^{\frac{\min\{1,\nu\}}{2}}}\right)^{k-2}\left(\frac{4\sqrt[k]{k}\sqrt{q^+}}{\sqrt{q^-}}\right)^k\notag\\
& \leqslant \frac{\left(2C_1^2 + 18\right)\left(\frac{16e^2q^+}{q^-}\right)}{q^-}\sum_{k=2}^{[\log N]}\left(\frac{2C_1 + 6}{\sqrt{q^{-}}N^{\frac{\min\{1,\nu\}}{2}}}\frac{4\sqrt[k]{k}\sqrt{q^+}}{\sqrt{q^-}}\right)^{k-2}\notag\\
&\leqslant\frac{32e^2q^+\left(C_1^2 + 9\right)}{\bigl(q^-\bigr)^2}\sum_{k=0}^{\infty}\left(\frac{2C_1 + 6}{q^{-}N^{\frac{\min\{1,\nu\}}{2}}}\frac{4\sqrt[k]{k}\sqrt{q^+}}{\sqrt{q^-}}\right)^{k-2}\notag\\
&\leqslant \frac{32e^2q^+\left(C_1^2 + 9\right)}{\bigl(q^-\bigr)^2}\sum_{k=0}^{\infty}\frac{1}{2^k}=\frac{64e^2q^+\left(C_1^2 + 9\right)}{\bigl(q^-\bigr)^2}.\label{E4_bound}
\end{align}

Using \eqref{H/mu_bound} we can bound $E_{5,N}(t)$ as well:
\begin{align}
\abs{E_{5,N}(t)}&=\frac{N}{(\log N)^4}\frac{\tilde{\mathcal{G}}_N(t)\frac{(\log N)^{2}}{\sqrt{N}}}{\mu^{\star}_{N}(t)}\sprod{\vk e, H_N(t)\,\vk e}\leqslant \frac{\tilde{\mathcal{G}}_N(t)}{(\log N)^2}\frac{6\sqrt{N}}{\sqrt{Nq^-}}\leqslant \frac{6}{\sqrt{q^-}}\tilde{\mathcal{G}}_N(t).\label{E5_bound}
\end{align}
The next step is to combine \eqref{E3_bound}, that resulted from the first step of the proof, with the bounds \eqref{E4_bound} and \eqref{E5_bound}.
We thus obtain the following upper bound on the residual term appearing in \eqref{zeta_last}, on $\Omega_N$ and for sufficiently large $N$,
\begin{align}
\abs{E_{6,N}(t)}\leqslant 2\sqrt{q^+} + \frac{8\mathfrak{C}^2\sqrt{q^+}}{q^-} + \frac{64e^2q^+\left(C_1^2 + 9\right)}{\bigl(q^-\bigr)^2} +\frac{6}{\sqrt{q^-}}\tilde{\mathcal{G}}_N(t).\label{E6_bound}
\end{align}
Finally, using that $\bar{\mu}_N(t)$ is deterministic, 
\begin{align}
\E{\mu_{N}^{\star}(t)} &= \bar{\mu}_N(t) + \E{\zeta_N(t)}\notag\\
&=\bar{\mu}_N(t) + \E{E_{6,N}(t)}\pk{\Omega_N}\frac{(\log N)^4}{N} + \E{\zeta_N(t)\mathbb{I}_{\Omega^{\rm c}_N}},\label{Emu*_last}
\end{align}
Hence, combining \eqref{zeta_def}, \eqref{zeta_last} and \eqref{Emu*_last},  on $\Omega_N$ and for sufficiently large $N$, 
\begin{align*}
\mu_{N}^{\star}(t) - \E{\mu_{N}^{\star}(t)} &= \sprod{\vk e,H_N(t)\,\vk e} + \big(E_{6,N}(t)-\E{E_{6,N}(t)}\pk{\Omega_N}\big)\frac{(\log N)^4}{N}+\E{\zeta_N(t)\mathbb{I}_{\Omega^{\rm c}_N}}.
\end{align*}
The claim follows by applying \eqref{E6_bound}, \eqref{Omega_bound} and the following bound, which can be obtained combining the definition \eqref{zeta_def} of $\zeta_N(t)$ with \eqref{bar_mu_rough_bound} and  \eqref{H_full_bound}:
\begin{align*}
    \abs{\E{\zeta_N(t)\mathbb{I}_{\Omega^{\rm c}_N}}} &\leqslant \E{\abs{\zeta_N(t)}\mathbb{I}_{\Omega^{\rm c}_N}}\\
    &\leqslant \E{\left(2\sqrt{Nq^+} + \left\Vert H_N(t) \right\Vert\right)\mathbb{I}_{\Omega^{\rm c}_N}} \\
    &\leqslant \left(2\sqrt{Nq^+} + C_1\sqrt{N}\right)\pk{\mathbb{I}_{\Omega^{\rm c}_N}}\\
    &\leqslant (2\sqrt{q^+} + C_1)e^{-\nu(\log N)^2 + \frac{1}{2}\log N}.
\end{align*}

\QED

\section{Proof of Theorem \ref{thm_expectation}}

Using \eqref{mu_from_mu*}, \netheo{thm_expectation} is a straightforward corollary of the following lemma.
\begin{lem}\label{expectation_full_expression}
For all $t\in[0,T]$
\begin{align*}
    \E{\mu^{\star}_{N}(t)}=\sqrt{Nq(t)}+\frac{1}{\sqrt{Nq(t)}} + \frac{\overline{\mathcal{G}}_N(t)}{\sqrt{p(t)(1-p(t))}}\left(\frac{1}{N\sqrt{N}}\right),
\end{align*}
where $\limsup_{N\to\infty}\sup_{t\in[0,T]}\abs{\overline{\mathcal{G}}_N(t)}$ is finite.
\end{lem}

The above lemma gives an uniform expansion of the expectation of the principal eigenvalue upto second order.

\medskip

\prooflem{expectation_full_expression} We follow the proof of \cite[Eqn (6.7)]{10.1214/11-AOP734}, but as the goal is to obtain a uniformly bounded function $\overline{\mathcal{G}}_N(t)$, and again the challenge lies in controlling the various terms of the series which have a lower order contribution.

This (lengthy) proof consists of three steps. 
As a first step, we quantify the contribution to $\mu_N^\star(t)$ corresponding to the event $\Omega_N$ defined in \nelem{uniform_event}, which eventually leads to the expansion \eqref{mu_expansion_5_for_expectation}.  
As a second step, we identify the candidate error term $\overline{\mathcal G}_N(t)$, now also including the contribution due to the complementary event $\Omega_{N}^{\rm c}$. In the third step we prove that $\sup_{t\in[0,T]}\abs{\overline{\mathcal{G}}(t,N)}$ is finite for sufficiently large $N$ . 

\medskip

\noindent {\sc Step I.}
By the representation \eqref{mu_recurrent_expansion}, 
\begin{align}
   {\mathfrak{m}}&=\sqrt{Nq(t)} + \frac{\sqrt{Nq(t)}}{\mathfrak{m} }\hfr{1}+\frac{\sqrt{Nq(t)}}{\mfr{2}}\hfr{2}+\frac{\sqrt{Nq(t)}}{\mfr{3}}\hfr{3}+\frac{1}{N\sqrt{N}}K_{1,N}(t).\label{mu_expansion_1_for_expectation}
\end{align}
where (for compactness) we locally write $\hfr{k}\equiv \mathfrak{h}_{k,N}(t):= \sprod{\vk e,H_N^k(t)\,\vk e}$, $\mathfrak{m}\equiv \mu_N^\star(t)$,  and
\begin{align*}
    K_{1,N}(t) &:= N^2\sqrt{q(t)}\sum_{k=4}^{\infty}\frac{\mathfrak{h}_{k}}{\mfr{k}}.
\end{align*}
The last term in the right hand side of \eqref{mu_expansion_1_for_expectation} can be dealt with in a straightforward manner. Indeed, using \eqref{H/mu_bound}, $|K_{1,N}(t)|$ can be bounded by a constant: for sufficiently large $N$,
\begin{align}
    \abs{K_1(t,N)} &\leqslant N^2\sqrt{q^+}\sum_{k=4}^{\infty}\left(\frac{\left\Vert  H_N(t) \right\Vert}{\mathfrak{m}}\right)^k\notag\\
    &\leqslant N^2\sqrt{q^+}\sum_{k=4}^{\infty}\left(\frac{6}{\sqrt{Nq^-}}\right)^k= \frac{6^4 N^2\sqrt{q^+}}{N^2(q^-)^2}\sum_{k=0}^{\infty}\left(\frac{6}{\sqrt{Nq^-}}\right)^k\notag\\
    &\leqslant \frac{1296 \sqrt{q^+}}{(q^-)^2}\sum_{k=0}^{\infty}\left(\frac{1}{2}\right)^k= \frac{2592 \sqrt{q^+}}{(q^-)^2}\label{K1_bound_2}.
\end{align}
To analyze the second, third, and fourth term in the right hand side of \eqref{mu_expansion_1_for_expectation}, we apply the following procedure. We write each of the fractions $\sqrt{Nq(t)}/\mathfrak{m}^i$ (where $i\in\{1,2,3\}$),  in terms of a quantity proportional to $1/\mathfrak{m}^{i-1}$ and quantities of the type $\sqrt{Nq(t)}/\mathfrak{m}^{j}$ with $j>i$. This is helpful because $\mathfrak{m}$ roughly behaves as $\sqrt{N}$, implying that all terms that are essentially behaving as $\sqrt{Nq(t)}/\mathfrak{m}^j$ for some $j>3$ are small relative to $1/N\sqrt{N}$. 
To do so, we add and subtract a term $\mathfrak{m}$ in the numerator, so that $1/\mathfrak{m}^{i-1}$ appears automatically, and for the difference $\mathfrak{m} - \sqrt{Nq(t)}$ we use the representation \eqref{mu_recurrent_expansion}. As a result, the first term cancels, and in what remains there is an additional $\mathfrak{m}$ in the denominator, thus achieving our goal. 

We now perform this procedure for the fraction in the second term in the right hand side of \eqref{mu_expansion_1_for_expectation}. We obtain, on $\Omega_N$, again as an application of the representation \eqref{mu_recurrent_expansion}: 
\begin{align}
    \frac{\sqrt{Nq(t)}}{\mathfrak{m} }&=1 - \frac{\mathfrak{m}  -\sqrt{Nq(t)}}{\mathfrak{m} }=1 - \frac{\sqrt{Nq(t)}}{\mathfrak{m} }\sum_{k=1}^{\infty}\frac{\hfr{k}}{\mfr{k}}
    \notag\\
    &=1 - \frac{\sqrt{Nq(t)}}{\mfr{2}}\hfr{1} - \frac{\sqrt{Nq(t)}}{\mfr{3}}\hfr{2} + \frac{1}{N\sqrt{N}}K_{2,N}(t).\label{second_term}
\end{align}
where
\begin{align*}
    K_{2,N}(t) := - \frac{N^2\sqrt{q(t)}}{\mathfrak{m} }\sum_{k=3}^{\infty}\frac{\hfr{k}}{\mfr{k}}.
\end{align*}
We can majorize $|K_{2,N}(t)|$ by a constant, similarly to how $|K_{1,N}(t)|$ was bounded. Indeed, for sufficiently large $N$ , using  \eqref{mu_rough_bound} and \eqref{H/mu_bound}, we find the upper bound 
\begin{align}
    \abs{K_{2,N}(t)} \leqslant  \frac{2N^2\sqrt{q^+}}{\sqrt{Nq^-}}\sum_{k=3}^{\infty}\left(\frac{6}{\sqrt{Nq^-}}\right)^k \leqslant  \frac{432\sqrt{q^+}}{\bigl(q^-\bigr)^2}\sum_{i=0}^{\infty}\frac{1}{2^k}=  \frac{864\sqrt{q^+}}{\bigl(q^-\bigr)^2}.\label{K2_bound_2}
\end{align}
Inserting \eqref{second_term} into \eqref{mu_expansion_1_for_expectation}, we thus obtain that,  on $\Omega_N$, 
\begin{align}
    \mathfrak{m}  &=\sqrt{Nq(t)} + \hfr{1} +\frac{\sqrt{Nq(t)}}{\mfr{2}}\bigl(\hfr{2} - \bigl(\hfr{1}\bigr)^2 \bigr)\label{mu_expansion_2_for_expectation}
    +\frac{\sqrt{Nq(t)}}{\mfr{3}}\bigl(\hfr{3} - \hfr{1}\hfr{2}\bigr)+\frac{1}{N\sqrt{N}}\bigl(K_{1,N}(t) + K_{2,N}(t)\,\hfr{1}\bigr).
\end{align}
Following the same procedure, we obtain for the fraction in the third term of the right hand side of \eqref{mu_expansion_2_for_expectation},
\begin{align}
    \frac{\sqrt{Nq(t)}}{\mfr{2}} &= \frac{\mathfrak{m} }{\mfr{2}} - \frac{\mathfrak{m} -\sqrt{Nq(t)}}{\mfr{2}}= \frac{1}{\mathfrak{m} } -  \frac{\sqrt{Nq(t)}}{\mfr{2}}\sum_{k=1}^{\infty}\frac{\hfr{k}}{\mfr{k}}=\frac{1}{\mathfrak{m}} -  \frac{\sqrt{Nq(t)}}{\mfr{3}}\hfr{1} + \frac{1}{N\sqrt{N}}K_{3,N}(t)\label{third_term}
\end{align}
where
\begin{align*}
    K_{3,N}(t) := -\frac{N^2\sqrt{q(t)}}{\mfr{2}}\sum_{k=2}^{\infty}\frac{\hfr{k}}{\mfr{k}}.
\end{align*}
Similarly to how we found the upper bound on $\abs{K_{1,N}(t)}$ and $\abs{K_{2,N}(t)}$, we now have
\begin{align}
 \abs{K_{3,N}(t)}\leqslant \frac{288\sqrt{q^+}}{\bigl(q^-\bigr)^2}\label{K3_bound_2}
\end{align}
Hence, upon combining \eqref{mu_expansion_2_for_expectation} and \eqref{third_term}, on $\Omega_N$,
\begin{align}
    \mathfrak{m} &=\sqrt{Nq(t)} + \hfr{1} +\frac{1}{\mathfrak{m}}\bigl(\hfr{2} - \bigl(\hfr{1}\bigr)^2 \bigr)\label{mu_expansion_3_for_expectation} +\frac{\sqrt{Nq(t)}}{\mfr{3}}\biggl(\hfr{3} - 2\hfr{1}\hfr{2} + \bigl(\hfr{1})^3\biggr)\\
    &\,\quad +\frac{1}{N\sqrt{N}}\left(K_{1,N}(t) + K_{2,N}(t)\,\hfr{1}+K_{3,N}(t)\bigl(\hfr{2} - \bigl(\hfr{1}\bigr)^2 \bigr)\right).\notag
\end{align}
Here, a term inversely proportional to ${\mathfrak{m}}$ appears, which we take care of later. Finally, for the fraction in the fourth term of the right hand side of \eqref{mu_expansion_3_for_expectation} we have
\begin{align}
    \frac{\sqrt{Nq(t)}}{\mfr{3}} &= \frac{1}{\mfr{2}} - \frac{\mathfrak{m} - \sqrt{Nq(t)}}{\mfr{3}}=\frac{1}{\mfr{2}} - \frac{\sqrt{Nq(t)}}{\mfr{3}}\sum_{k=1}^{\infty}\frac{\hfr{k}}{\mfr{k}}=\frac{1}{\mfr{2}} +
    \frac{1}{N\sqrt{N}}K_{4,N}(t),\label{fourth_term}
\end{align}
where
\begin{align*}
    K_{4,N}(t) := - \frac{\sqrt{Nq(t)}}{\mfr{3}}\sum_{k=1}^{\infty}\frac{\hfr{k}}{\mfr{k}},
\end{align*}
which satisfies
\begin{align}
    \abs{K_{4,N}(t)} \leqslant  \frac{96\sqrt{q^+}}{\bigl(q^-\bigr)^2}\label{K4_bound_2}.
\end{align}
Inserting \eqref{fourth_term} into \eqref{mu_expansion_3_for_expectation}, we have now arrived at
\begin{align}
    \mathfrak{m} &=\sqrt{Nq(t)} + \hfr{1} +\frac{1}{\mathfrak{m}}\bigl(\hfr{2} - \bigl(\hfr{1}\bigr)^2 \bigr)\label{mu_expansion_4_for_expectation} +\frac{1}{\mfr{2}}\biggl(\hfr{3} - 2\hfr{1}\hfr{2} + \bigl(\hfr{1})^3\biggr) +\frac{1}{N\sqrt{N}}K_{5,N}(t),
\end{align}
where
\begin{align*}
    K_{5,N}(t) &:= K_{1,N}(t) + K_{2,N}(t)\,\hfr{1} +K_{3,N}(t)\bigl(\hfr{2} - \bigl(\hfr{1}\bigr)^2 \bigr)+ K_{4,N}(t)\biggl(\hfr{3} - 2\hfr{1}\hfr{2} + \bigl(\hfr{1})^3\biggr) .
\end{align*}
The next step is to combine \eqref{A1} with the obtained bounds \eqref{K1_bound_2}, \eqref{K2_bound_2}, \eqref{K3_bound_2}, and  \eqref{K4_bound_2}, in order to obtain the following uniform bound: on $\Omega_N$,
\begin{align}
    \abs{K_{5,N}(t)} \leqslant \frac{2592 \sqrt{q^+}}{(q^-)^2} + 3\frac{864\sqrt{q^+}}{\bigl(q^-\bigr)^2} +18\frac{288\sqrt{q^+}}{\bigl(q^-\bigr)^2} + 108 \frac{96\sqrt{q^+}}{\bigl(q^-\bigr)^2}= \frac{20736 \sqrt{q^+}}{(q^-)^2}.\label{K5_bound_2}
\end{align}
We proceed by analyzing the terms in \eqref{mu_expansion_4_for_expectation} which are proportional to $1/{\mathfrak{m}}$ and $1/\mfr{2}$.  Recalling that by \nelem{mu_full_expression} we have that $1/{\mathfrak{m}}\approx 1/\sqrt{Nq(t)}$ and $1/\mfr{2}\approx 1/Nq(t)$,  using the expansion given in \nelem{lem_mu_recurrent_expansion},
\begin{align}
    \frac{1}{\mathfrak{m}} &= \frac{1}{\sqrt{Nq(t)}} - \frac{\mathfrak{m} - \sqrt{Nq(t)}}{\sqrt{Nq(t)}\,\mathfrak{m}}=\frac{1}{\sqrt{Nq(t)}} - \frac{1}{\mathfrak{m}}\sum_{k=1}^{\infty}\frac{\hfr{k}}{\mfr{k}}=\frac{1}{\sqrt{Nq(t)}} - \frac{\hfr{1}}{\mfr{2}} + \frac{1}{N\sqrt{N}}K_{6,N}(t)\label{fifth_term},
    \end{align}
    and likewise
    \begin{align}
    \frac{1}{\mfr{2}} & = \frac{1}{Nq(t)} - \left(\frac{1}{\sqrt{Nq(t)}} + \frac{1}{\mathfrak{m}}\right)\left(\frac{1}{\sqrt{Nq(t)}} - \frac{1}{\mathfrak{m}}\right)\notag\\
    & = \frac{1}{Nq(t)} - \left(\frac{1}{\sqrt{Nq(t)}} + \frac{1}{\mathfrak{m}}\right)\left(\frac{\hfr{1}}{\mfr{2}} - \frac{1}{N\sqrt{N}}K_{6,N}(t)\right)= \frac{1}{Nq(t)} + \frac{1}{N\sqrt{N}}K_{7,N}(t)\label{sixth_term}
\end{align}
for the functions
\begin{align*}
    K_{6,N}(t) &:= - \frac{N\sqrt{N}}{\mathfrak{m}}\sum_{k=2}^{\infty}\frac{\hfr{k}}{\mfr{k}},\quad\quad
    K_{7,N}(t) := - \left(\frac{\sqrt{N}}{\sqrt{Nq(t)}} + \frac{\sqrt{N}}{\mathfrak{m}}\right)\left(\frac{N\,\hfr{k}}{\mfr{2}} - \frac{N}{N\sqrt{N}}K_{6,N}(t)\right).
\end{align*}
These functions  can be uniformly bounded in the (by now) standard manner:
\begin{align}
    \abs{K_{6,N}(t)} &\leqslant \frac{2N\sqrt{N}}{\sqrt{Nq^-}}\sum_{k=2}^{\infty}\left(\frac{6}{\sqrt{Nq^-}}\right)^k\leqslant \frac{72N\sqrt{N}}{Nq^-\sqrt{Nq^-}}\sum_{k=0}^{\infty}\frac{1}{2^k}= \frac{144}{q^-\sqrt{q^-}},\label{K6_bound_2}\\
    \abs{K_{7,N}(t)} &\leqslant \left(\frac{1}{\sqrt{q^-}} + \frac{2}{\sqrt{q^-}}\right)\left(\frac{12}{q^-} + \frac{144}{q^-\sqrt{q^-}}\right) = \frac{36\sqrt{q^-} + 432}{\bigl(q^-\bigr)^2}\label{K7_bound_2}
\end{align}
Hence, combining \eqref{mu_expansion_4_for_expectation} with \eqref{fifth_term} and \eqref{sixth_term}, we have obtained our final expansion: on $\Omega_N$,
\begin{align}
    {\mathfrak{m}} &=\sqrt{Nq(t)} + \hfr{1} + \frac{1}{\sqrt{Nq(t)}}\bigl(\hfr{2} - \bigl(\hfr{1}\bigr)^2 \bigr)\label{mu_expansion_5_for_expectation} +\frac{1}{Nq(t)}\biggl(\hfr{3} - 3\hfr{1}\hfr{2} + 2\bigl(\hfr{1}\bigr)^3\biggr)+\frac{1}{N\sqrt{N}}K_{8,N}(t),
\end{align}
where
\begin{align*}
    K_{8,N}(t) &:= K_{5,N}(t) + K_{6,N}(t)\bigl(\hfr{2} - \bigl(\hfr{1}\bigr)^2 \bigr)+ K_{7,N}(t)\biggl(\hfr{3} - 3\hfr{1}\hfr{2} + 2\bigl(\hfr{1}\bigr)^3\biggr).
\end{align*}
We then combine \eqref{K5_bound_2}, \eqref{K6_bound_2} and \eqref{K7_bound_2} with \eqref{A1}. We thus obtain, for any $t\in[0,T]$ and large enough $N$, the following bound: on $\Omega_N$,
\begin{align}
    \abs{K_{8,N}(t)}&\leqslant \frac{20736 \sqrt{q^+}}{(q^-)^2} + 18\frac{144}{q^-\sqrt{q^-}} + (27+81+54)\frac{36\sqrt{q^-} + 432}{\bigl(q^-\bigr)^2}\notag\\
    &= \frac{20736 \sqrt{q^+} + 8424\sqrt{q^-} + 69984}{(q^-)^2}.\label{K8_bound_2}
\end{align}
{\sc Step II.} We proceed by using \eqref{mu_expansion_5_for_expectation} to analyze our target quantity $\E{\mu_N^\star(t)}$, by also taking into account the contribution due to $\Omega_N^{\rm c}$.
To this end, we first define the objects $\mathfrak{h}_{k,N}^\circ(t):=\mathfrak{h}_{k,N}(t) \,{\mathbb I}\{\Omega_N\}$ and $\mathfrak{h}_{k,N}^{\rm c}(t):=\mathfrak{h}_{k,N}(t) \,{\mathbb I}\{\Omega_N^{\rm c}\}$.
From \eqref{S2} we know that 
\begin{align*}
 \E{\mathfrak{h}_{1,N}(t)}&=
    \E{\sprod{\vk e,H_N(t)\,\vk e}} = 0,\\
    \E{\mathfrak{h}_{2,N}(t)}&=\E{\sprod{\vk e,H^2_N(t)\,\vk e}} = 1.
    \end{align*}
Hence, combining this with \eqref{mu_expansion_5_for_expectation}, we obtain
\begin{align*}
    \E{\mu_N^\star(t)} &=\E{\mu_N^\star(t)\mathbb{I}\{{\Omega_N}\}} + \E{\mu_N^\star(t)\mathbb{I}\{{\Omega_N^{\rm c}}\}}\\
    &=\sqrt{Nq(t)}\,\pk{\Omega_N} - \E{\mathfrak{h}_{1,N}^{\rm c}(t)} \notag\\
    &\quad+\, \frac{\pk{\Omega_N}-\E{\mathfrak{h}_{2,N}^{\rm c}(t)}-\E{\big(\mathfrak{h}_{1,N}^\circ(t)\big)^2}}{\sqrt{Nq(t)}}\notag\\
&\quad+\frac{\E{\mathfrak{h}_{3,N}^\circ(t)}-3\E{\mathfrak{h}_{1,N}^\circ(t)\,\mathfrak{h}_{2,N}^\circ(t)}+2\E{\big(\mathfrak{h}_{1,N}^\circ(t)\big)^2}}{Nq(t)}\\
&\quad+\, \frac{\E{K_{8,N}(t)\,{\mathbb I}\{\Omega_N\}}}{N\sqrt{N}}+ \E{\mu_N^\star(t)\mathbb{I}\{{\Omega_N^{\rm c}}\}}\\
   &=\sqrt{Nq(t)}+\frac{1}{\sqrt{Nq(t)}} + \frac{\overline{\mathcal{G}}_N(t)}{\sqrt{p(t)(1-p(t))}}\left(\frac{1}{N\sqrt{N}}\right)\notag,
\end{align*}
where the function $\overline{\mathcal{G}}_N(t)$ is the sum of nine terms: with $\bar p(t):=\sqrt{p(t)(1-p(t))}$,
\begin{align}
    \overline{\mathcal{G}}_N(t) &:= -\bar p(t)\left(N^2\sqrt{q(t)}+\frac{N}{\sqrt{q(t)}}\right)\pk{\Omega_N^{\rm c}}\notag-N\sqrt{N} \,\bar p(t) \,\E{\mathfrak{h}_{1,N}^{\rm c}(t)}\\
    \notag &\quad -N\,\bar p(t) \,\frac{\E{\mathfrak{h}_{2,N}^{\rm c}(t)}}{\sqrt{q(t)}}-N\,\bar p(t)\,\frac{\E{\big(\mathfrak{h}_{1,N}^{\circ}(t)\big)^2}}{\sqrt{q(t)}} \\
    \notag &\quad 
    +\,\sqrt{N}\,\bar p(t)\frac{\E{\mathfrak{h}_{3,N}^\circ(t)}}{q(t)} - 3\sqrt{N}\,\bar p(t)\frac{ \E{\mathfrak{h}_{1,N}^\circ(t)\,\mathfrak{h}_{2,N}^\circ(t) } }{q(t)} +2 \sqrt{N}\,\bar p(t)\frac{ \E{\big(\mathfrak{h}_{1,N}^\circ(t)\big)^2  }}{q(t)}\\
    &\quad +\, \bar p(t) \,\E{K_{8,N}(t)\,{\mathbb I}\{\Omega_N\}}+N\sqrt{N}\,\bar p(t) \,\E{\mu_N^\star(t)\mathbb{I}\{{\Omega_N^{\rm c}}\}}.
    \label{G_bar}
\end{align}
{\sc Step III.} We have thus arrived at the last part of the proof, where we are to bound the above sum, uniformly in $t\in[0,T]$ for sufficiently large $N$. We do this by consider each of the nine terms in the right hand side of \eqref{G_bar} separately. 
The first seven terms of \eqref{G_bar} can be bounded using \eqref{S2} and \eqref{Omega_bound} for sufficiently large $N$ , whereas the last two terms have to be dealt with differently. We throughout use the notation $p^\circ:=\sqrt{p^+ (1-p^-)}$.

\medskip 

\noindent {\it First term.}
This term can be bounded using \eqref{Omega_bound}: for $N$ sufficiently large,
\begin{align*}
    \abs{\,\bar p(t)\left(N^2\sqrt{q(t)}+\frac{N}{\sqrt{q(t)}}\right)\pk{\Omega_N^{\rm c}}}&\leqslant p^\circ \left(N^2\sqrt{q^+}+\frac{N}{\sqrt{q^-}}\right)e^{-\nu(\log N)^2}\leqslant 1.
\end{align*}
In the next six terms we rely on the following bounds, each of them following from \eqref{S2}:
\begin{align}
    \abs{\mathfrak{h}_{1,N}(t)}&\leqslant C_1\sqrt{N},&&\label{sprod_bound_1}\\
    \abs{{\mathfrak{h}_{2,N}(t)}}&\leqslant C_1^2 N,&&\label{sprod_bound_2}\\
    \abs{(\mathfrak{h}_{1,N}(t))^2}&\leqslant C_1^2N,\qquad &\abs{\E{(\mathfrak{h}_{1,N}(t))^2} }&\leqslant {3}/{N},\label{sprod_bound_3}\\
    \abs{{\mathfrak{h}_{3,N}(t)}}&\leqslant C_1^3 N\sqrt{N},\qquad &\abs{\E{\mathfrak{h}_{3,N}(t)}} &\leqslant {8C_2^3}/{\sqrt{N}},\label{sprod_bound_4}\\
    \abs{\mathfrak{h}_{1,N}(t)\mathfrak{h}_{2,N}(t)}&\leqslant C_1^3 N\sqrt{N} \qquad &\abs{\E{\mathfrak{h}_{1,N}(t)\mathfrak{h}_{2,N}(t)}}&\leqslant {8C_2^3}/{(N\sqrt{N})},\label{sprod_bound_5}\\
    \abs{(\mathfrak{h}_{1,N}(t))^3}&\leqslant C_1^2 N\sqrt{N}\qquad &\abs{\E{(\mathfrak{h}_{1,N}(t))^3}} &\leqslant  {8C_2^3}/{(N^2\sqrt{N})}.\label{sprod_bound_6}
\end{align}

The left column presented above straightforwardly follows from the first inequality of \eqref{S2} after simply expand the inner product in the definition of $\mathfrak{h}_k(t)$. To show the inequalities in the right column, we  use the fact that the entries of $H_N(t)$ are independent and centered, meaning that if we expand the products under the expectation, a considerable number of terms equals zero. In the right inequality of \eqref{sprod_bound_3} we obtain 
\begin{align*}
    \abs{\E{(\mathfrak{h}_{1,N}(t))^2}} \leqslant \frac{1}{N}\sum_{i,j=1}^{N^2}\left( 2\abs{\E{h_{i,j}^2(t)}} + 2\abs{\E{h_{i,j}h_{j,i}}}\right) =\frac{3}{N},
\end{align*}
where the last inequality comes from the third equation of \eqref{S2}.
By the same token, in the right inequality of \eqref{sprod_bound_4} we have
\begin{align*}
    \abs{\E{\mathfrak{h}_{3,N}(t)}} \leqslant \frac{1}{N}\sum_{i,j=1}^{N^2} \abs{8\E{h_{i,j}^3(t)}} \leqslant\frac{8C_2^3}{\sqrt{N}},
\end{align*}
in the right inequality of \eqref{sprod_bound_5} we have
\begin{align*}
    \abs{\E{\mathfrak{h}_{1,N}(t)\mathfrak{h}_{2,N}(t)}} \leqslant \frac{1}{N^2}\sum_{i,j=1}^{N^2} \abs{8\E{h_{i,j}^3(t)}} \leqslant\frac{8C_2^3}{N\sqrt{N}},
\end{align*}
and in the right inequality of \eqref{sprod_bound_6} we have
\begin{align*}
    \abs{\E{(\mathfrak{h}_{1,N}(t))^3}} \leqslant \frac{1}{N^3}\sum_{i,j=1}^{N^2} \abs{8\E{h_{i,j}^3(t)}} \leqslant\frac{8C_2^3}{N\sqrt{N}},
\end{align*}
where in this three cases the last inequalities come from the last inequality of \eqref{S2} for $r=3$.

\noindent
{\it Second term.} We combine \eqref{Omega_bound} with \eqref{sprod_bound_1}: for $N$ sufficiently large,
\begin{align*}
    \abs{N\sqrt{N} \,\bar p(t) \,\E{\mathfrak{h}_{1,N}^{\rm c}(t)}} &\leqslant N\sqrt{N}p^\circ \E{C_1\sqrt{N}\,\mathbb{I}\{{\Omega^{\rm c}_N}\}}\leqslant C_1N^2p^\circ e^{-\nu(\log N)^2}\leqslant 1.
\end{align*}
{\it Third term.} We combine \eqref{Omega_bound} with \eqref{sprod_bound_2}: for $N$ sufficiently large,
\begin{align*}
    \abs{N\,\bar p(t) \,\frac{\E{\mathfrak{h}_{2,N}^{\rm c}(t)}}{\sqrt{q(t)}}}&\leqslant Np^\circ \,\frac{\E{C_1^2 N\,\mathbb{I}\{{\Omega_N^{\rm c}}\}}}{\sqrt{q^-}}\leqslant C_1^2N^2p^\circ \frac{e^{-\nu(\log N)^2}}{\sqrt{q^-}}\leqslant 1.
\end{align*}
{\it Fourth  term.} We combine  \eqref{Omega_bound} and \eqref{sprod_bound_3}: for $N$ sufficiently large,
\begin{align*}
    \abs{N\,\bar p(t)\,\frac{\E{\big(\mathfrak{h}_{1,N}^{\circ}(t)\big)^2}}{\sqrt{q(t)}}}&\leqslant Np^\circ \frac{\E{\big(\mathfrak{h}_{1,N}(t)\big)^2}}{\sqrt{q^-}}+ Np^\circ \frac{\E{\big(\mathfrak{h}_{1,N}^{\rm c}(t)\big)^2}}{\sqrt{q^-}}\\
    &\leqslant Np^\circ \frac{4}{N\sqrt{q^-}}+ Np^\circ \frac{C_1^2Ne^{-\nu(\log N)^2}}{\sqrt{q^-}}\leqslant \frac{3p^\circ }{\sqrt{q^-}} + 1.
\end{align*}
{\it Fifth term.} We combine \eqref{Omega_bound} and \eqref{sprod_bound_4}:
for $N$ sufficiently large,
\begin{align*}
    \abs{\sqrt{N}\,\bar p(t)\frac{\E{\mathfrak{h}_{3,N}^\circ(t)}}{q(t)}}&\leqslant \sqrt{N}p^\circ \frac{\abs{\E{\mathfrak{h}_{3,N}(t)}}}{q^-}+\sqrt{N}p^\circ \frac{\abs{\E{\mathfrak{h}^{\rm c}_{3,N}(t)}}}{q^-}\\
    &\leqslant \sqrt{N}p^\circ \frac{8C_2^3}{\sqrt{N}q^-}+p^\circ \frac{C_1^3N^2e^{-\nu(\log N)^2}}{q^-}\leqslant p^\circ \frac{8C_2^3}{q^-} + 1.
\end{align*}
{\it Sixth term.} We combine \eqref{Omega_bound} and \eqref{sprod_bound_5}:
for $N$ sufficiently large,
\begin{align*}
    \abs{\sqrt{N}\,\bar p(t)\frac{\E{\mathfrak{h}_{1,N}^\circ(t)\,\mathfrak{h}_{2,N}^\circ(t)  }}{q(t)}}&\leqslant \sqrt{N}p^\circ \frac{\E{\mathfrak{h}_{1,N}(t)\,\mathfrak{h}_{2,N}(t)  }}{q^-}+\sqrt{N}p^\circ \frac{\E{\mathfrak{h}_{1,N}^{\rm c}(t)\,\mathfrak{h}_{2,N}^{\rm c}(t)  }}{q^-}\\
    &\qquad\leqslant \sqrt{N}p^\circ \frac{8C_2^3}{N\sqrt{N}q^-}+\sqrt{N}p^\circ \frac{C_1^3 N\sqrt{N}e^{-\nu(\log N)^2}}{q^-}\leqslant 2.
\end{align*}
{\it Seventh term.} We combine \eqref{Omega_bound} and \eqref{sprod_bound_6}:
for $N$ sufficiently large,
\begin{align*}
    \sqrt{N}\,\bar p(t)\frac{ \E{\big(\mathfrak{h}_{1,N}^\circ(t)\big)^2  }}{q(t)}&\leqslant \sqrt{N}p^\circ \frac{\E{\big(\mathfrak{h}_{1,N}(t)\big)^2  }}{q^-}+\sqrt{N}p^\circ \frac{\E{\big(\mathfrak{h}_{1,N}^{\rm c}(t)\big)^2  }}{q^-}\\
    &\leqslant \sqrt{N}p^\circ \frac{8C_2^3}{N^2\sqrt{N}q^-}+\sqrt{N}p^\circ \frac{C_1^2 N\sqrt{N}e^{-\nu(\log N)^2}}{q^-}\leqslant 2.
\end{align*}
{\it Eighth term.} Using \eqref{K8_bound_2}, for $N$ sufficiently large,
\begin{align*}
    \bar p(t) \,\E{K_{8,N}(t)\,{\mathbb I}\{\Omega_N\}}&\leqslant p^\circ \,\E{\abs{K_{8,N}(t)}\,{\mathbb I}\{\Omega_N\}}\leqslant p^\circ \frac{20736 \sqrt{q^+} + 8424\sqrt{q^-} + 69984}{(q^-)^2}.
\end{align*}
{\it Ninth term.} We combine \eqref{Omega_bound} with the definitions of $\mu_N^{\star}(t)$ and $A_{N}^{\star}(t)$ that were given in \eqref{A*_def}: for $N$ sufficiently large,
\begin{align*}
    N\sqrt{N}\,\bar p(t) \,\E{\mu_N^\star(t)\mathbb{I}\{{\Omega_N^{\rm c}}\}}&\leqslant N\sqrt{N}\,\bar p(t)\,\E{\left\Vert A_N^{\star}(t)\right\Vert\mathbb{I}\{{\Omega_N^{\rm c}}\}}\\
    &= N\,\E{\left\Vert A_N(t)\right\Vert\mathbb{I}\{{\Omega}_N^{\rm c}\}}\leqslant N\,\E{N\mathbb{I}\{{\Omega_N^{\rm c}}\}}\leqslant N^2e^{-\nu(\log N)^2}\leqslant 1.
\end{align*}
Combining all the bounds above we obtain that there is a constant upper bound on $\overline{\mathcal{G}}_N(t)$ that applies to any $t\in [0,T]$
and $N$ sufficiently large. Hence, the claim follows.
\QED

\section{Discussion and concluding remarks}
In this paper we have considered the dynamic version of the conventional Erd\H{o}s-R\'enyi model, in which the edges independently alternate between being absent and present, with on- and off-times that are exponentially distributed. The main result is a functional central limit theorem.

\medskip

Extensions can be explored further in many directions. We discuss five of them.
\begin{itemize}
    \item[$\circ$] We saw that for our model the principal eigenvalue effectively inherits the per-edge correlation structure (see Remark \ref{R1}). One may wonder whether such a principle holds for general on- and off-time distributions. Observe that this would require a rather different tightness proof, as the one we have set up in this paper extensively uses the underlying Markovian structure. 
    \item[$\circ$] A branch of the dynamic random graph literature has a focus on estimating (the parameters pertaining to) the underlying random dynamics from partial information; see e.g.\ \cite{mandjeswang2024} and references therein. In the context of the present paper a concrete question would be: observing the principal eigenvalue (say) at equidistant points in time, can we estimate the parameters $\LON$ and $\LOFF$?  Given the estimates we have of $\E{\mu_N(t)\,\mu_N(t+\delta)}$, one would anticipate that a method of moments can be developed. 
    \item[$\circ$] In many real-life examples, edge dynamics is affected by a common background process. This mechanism, referred to in the literature as {\it regime switching}, has a wide application potential, as it can model the impact of weather conditions, temperature, or any other exogenous process. Earlier results \cite{mandjes2019dynamic} suggest that when the background process is of Markovian type and when time is scaled appropriately, the principal eigenvalue's limiting process (in the central-limit regime, that is) remains effectively of the same Ornstein-Uhlenbeck nature. 
    \item[$\circ$] In this paper we studied the principal eigenvalue process for the dynamic version of the Erd\H{o}s-R\'enyi random graph. It is natural to consider dynamic versions of other standard random graph models as well, where a starting point could be the Chung-Lu model; cf.\ \cite{zhang2017random, dionigi2023central}.
    \item[$\circ$] A last topic for follow up research could focus on functional central limit theorems for other graph-related quantities. In this context one could in particular think of subgraph counts, such as the number of wedges or triangles. Here one should somehow deal with the complication that such subgraph counts are not Markov (except for the number of edges). 
    
\end{itemize}

\bibliographystyle{imsart-number} 
\bibliography{reference.bib}

\appendix

\section{Auxiliary results}

\medskip

To show \nelem{uniform_event} we need Lemmas \ref{A1_in_distribution} and \ref{A2_in_distribution}, which are actually the uniform analog of \cite[Lemmas 4.3 and 6.5]{10.1214/11-AOP734}. We shall show that under \neprop{model_basic_properties} the bounds obtained in \cite{10.1214/11-AOP734} are actually already uniform. 

\begin{lem}\label{A1_in_distribution}
    There exist a constant $\nu_1$ such that, for any $t\in[0,T]$ and sufficiently large $N$,
    \begin{align*}
        \pk{\left\Vert H_N(t)\right\Vert\leqslant 2+\frac{(\log N)^{2}}{N^{1/4}}}\geqslant 1- e ^{-\nu_1(\log N)^{2}}.
    \end{align*}    
\end{lem}

\prooflem{A1_in_distribution} Using \cite[(A.24)]{10.1214/11-AOP734} and the properties of the matrix $H_N(t)$, we obtain for $C^\star :=2\sqrt{C_1\sqrt{p^+}}$ and any even integer $k<C^\star (N)^{1/4}$
\begin{align*}
\E{\left(\left\Vert H_N(t)\right\Vert\right)^k}\leqslant \abs{\E{\operatorname{Tr}H_N^k(t)}}\leqslant 3N2^k.
\end{align*}
Hence, for any even $k\in\big(\frac{C^\star }{2}N^{1/4},\, C^\star N^{1/4}\big)$, using the Markov inequality,
\begin{align*}
    \pk{\left\Vert H_N(t)\right\Vert > 2+\frac{(\log N)^{2}}{N^{1/4}}} &= \pk{\left(\left\Vert H_N(t)\right\Vert\right)^k > \left(2+\frac{(\log N)^{2}}{N^{1/4}}\right)^k}\leqslant \frac{3N2^k}{\left(2+\frac{(\log N)^{2}}{N^{1/4}}\right)^k}\\
    &= \exp\left(\log(3N) + k\left(\log 2 - \log\left(2+\frac{(\log N)^2}{N^{1/4}}\right)\right)\right)\\
    &\leqslant \exp\left(\log(3N) - k\frac{(\log N)^2)}{3N^{1/4}}\right)\\
    &\leqslant \exp\left(\log(3N) - \frac{C^\star }{6}(\log N)^2)\right)\leqslant \exp\left( - \frac{C^\star }{12}(\log N)^2)\right).
\end{align*}
Hence, the claim follows for $\nu_1 := C^\star /12 = \sqrt{C_1\sqrt{p^+}}/6$.
\QED

\begin{lem}\label{A2_in_distribution}
   There exist positive constants $\nu_2$, $\mathfrak{C}$ such that for any $k\in\{1,\ldots, \lceil \log N\rceil\} $, any $t\in[0,T]$ and sufficiently large $N$,
    \begin{align*}
        \pk{\abs{\sprod{\vk e,H^k_N(t)\,\vk e}-\E{\sprod{\vk e, H^k_N(t)\,\vk e}}}<\mathfrak{C}^{k}\frac{(\log N )^{2k}}{\sqrt{N}}}\geqslant 1- e ^{-\nu_2(\log N)^{2}}.
    \end{align*}    
\end{lem}

\prooflem{A2_in_distribution}
Following the line of the proof of \cite[Lemma 6.5]{10.1214/11-AOP734} we obtain that, for any $t\in[0,T]$ and any $k,p\in\N$,
\begin{align*}
    \E{\abs{\sprod{\vk e, H_N^k(t)\vk e} - \E{\sprod{\vk e, H_N^k(t)\vk e}}}^p} \leqslant \left(\frac{(2Ckp)^k}{\sqrt{N}}\right)^p
\end{align*}
for the constant $C = \max(1,C_2),$
where $C_2$ is defined in \eqref{S2}.
Hence, using that $k\leqslant \log N + 1$, defining 
\begin{align*}
    \mathfrak{C} &:= 4eC,\:\:\:\:\:p := \left\lceil \frac{(\log N)^2}{k}\right\rceil,
\end{align*}
and applying the Markov inequality, we obtain that
\begin{align*}
    \pk{\abs{\sprod{\vk e,H^k_N(t)\,\vk e}-\E{\sprod{\vk e, H^k_N(t)\,\vk e}}}\geqslant\mathfrak{C}^{k}\frac{(\log N )^{2k}}{\sqrt{N}}}&\leqslant \frac{\left(\frac{(2Ckp)^k}{\sqrt{N}}\right)^p}{\left(\mathfrak{C}^{k}\frac{(\log N )^{2k}}{\sqrt{N}}\right)^p}\\
    & = \left(\frac{kp}{2e(\log N)^{2}}\right)^{pk} \leqslant e^{-(\log N)^2}.
\end{align*}
Hence, the claim follows for $\nu_2 = 1$.
\QED

\medskip

\prooflem{lem_mu_recurrent_expansion} 
Let $\vk v_N(t)$ is the eigenvector of the matrix $A^{\star}_N(t)$ corresponding to eigenvalue $\mu^{\star}_{N}(t)$, meaning that $A^{\star}_N(t)\, \vk v_N(t) = \mu^{\star}_{N}(t)\,\vk v_N(t).$
Using the representation \eqref{A_in_H}, we thus have that
\begin{align*}
    \left(H_N(t)+\sqrt{Nq(t)}\,\mathcal{E}_N\right)\vk v_N(t) = \mu^{\star}_{N}(t)\,\vk v_N(t),
\end{align*}
implying
\begin{align}
    (\mu^{\star}_{N}(t) - H_N(t))\,\vk v_N(t) = \sqrt{Nq(t)}\sprod{\vk e,\vk v_N(t)}\vk e.\label{tmp1}
\end{align}
Given that for sufficiently large $N$, $\mu^{\star}{N}(t) = \left\Vert A_N^*(t)\right\Vert > 4$, and from \eqref{A1}, for sufficiently large $N$, $\left\Vert H_N(t)\right\Vert < 3$, we observe that $\left\Vert\mu^{\star}{N}(t) - H_N(t)\right\Vert > 1$, which is indeed greater than 0. Hence, the matrix $\mu^{\star}_{N}(t) - H_N(t)$ is invertible,
so that for some function $\bar{Q}_N(t)$ we can write
\begin{align}
\vk v_N(t) = \bar{Q}_N(t)(\mu^{\star}_{N}(t) - H_N(t))^{-1}\vk e\label{tmp2}
\end{align}
Inserting \eqref{tmp2} into \eqref{tmp1}, we obtain
\begin{align*}
    1 = \sqrt{Nq(t)}\sprod{\vk e, (\mu^{\star}_{N}(t) - H_N(t))^{-1}\vk e}.
\end{align*}
Now \eqref{mu_recurrent_expansion} follows from the fact that
\begin{align*}
    \bigl(\mu^{\star}_{N}(t)-H_N(t)\bigr)^{-1} =\frac{1}{\mu^{\star}_{N}(t)}\sum_{k= 0}^\infty\left(\frac{H_N(t)}{\mu^{\star}_{N}(t)}\right)^k,
\end{align*}
which proves our claim.
\QED

\medskip

\prooflem{uniform_event}
Let $\nu_1$ be as defined in \nelem{A1_in_distribution} and $\nu_2$ be as defined in \nelem{A2_in_distribution}. Define $\nu^{\star}:=\min(\nu_1,\nu_2)$ and 
\begin{align*}
t_j = jT\,e^{-\frac{\nu^{\star}}{2}(\log N)^2}
\end{align*}
for $j\in\mathcal{J}:=\{0,\ldots, [e^{\frac{\nu^{\star}}{2} (\log N)^2}]\}$. Then, for sufficiently large $N$, using \nelem{A1_in_distribution},
\begin{align}
\pk{\exists j\in\mathcal{J}\colon \left\Vert H_N(t_j)\right\Vert > 2+\frac{(\log N)^2}{N^{1/4}}}&\leqslant \sum_{j=1}^{[e^{\frac{\nu^{\star}}{2} (\log N)^2}]}e ^{-\nu_2(\log N)^2}\notag\\
&\leqslant e ^{\frac{\nu^{\star}}{2}(\log N)^2}e ^{-\nu_2(\log N)^2}\leqslant e ^{-\frac{\nu^{\star}}{2}(\log N)^2}.\label{O1_bound}
\end{align}
By the same reasoning , using \nelem{A2_in_distribution}, for sufficiently large $N$,
\begin{align}
        \pk{\exists j\in\mathcal{J}\colon \abs{\sprod{\vk e,H^k_N(t_j)\,\vk e}-\E{\sprod{\vk e,H^k_N(t_j)\,\vk e}}}\geqslant \mathfrak{C}^{k}\frac{(\log N)^{2k}}{\sqrt{N}}}\leqslant e ^{-\frac{\nu^{\star}}{2}(\log N)^2}\label{O2_bound}
    \end{align}    
for the constants $\mathfrak{C}$ defined in \nelem{A2_in_distribution}.
Using \eqref{A0} from \neprop{model_basic_properties}, for sufficiently large  $N$,
\begin{align}
    \pk{\exists t^{\star}_1,t^{\star}_2\in\mathcal{T}_N: \abs{t^{\star}_1 - t^{\star}_2}\leqslant Te^{-\frac{\nu^{\star}}{2}(\log N )^2}}&\leqslant  CN^{4}Te^{-\frac{\nu^{\star}}{2}(\log N)^2}\notag\\
    &\leqslant e^{\log C + 5\log N -\frac{\nu^{\star}}{2}(\log N)^2}\leqslant e^{-\left(\frac{\nu^{\star}}{4}\right)(\log N)^2}.\label{O3_bound}
\end{align}
Define the events
\begin{align*}
    \Omega_{1,N}&:=\left\{\forall j\in\mathcal{J}\colon \left\Vert H_N(t_j)\right\Vert \leqslant 2+\frac{(\log N)^2}{N^{1/4}}\right\},\\
    \Omega_{2,N}&:=\left\{\forall j\in\mathcal{J}\colon \abs{\sprod{\vk e,H^k_N(t_j)\,\vk e}-\E{\sprod{\vk e,H^k_N(t_j)\,\vk e}}}< \mathfrak{C}^{k}\frac{(\log N )^{2k}}{\sqrt{N}}\right\},\\
    \Omega_{3,N}&:=\left\{\forall t^{\star}_1,t^{\star}_2\in\mathcal{T}_N: \abs{t^{\star}_1 - t^{\star}_2} >  Te^{-\frac{\nu^{\star}}{2}(\log N )^2}\right\}.
\end{align*}
Then for $\Omega_N:=\Omega_{1,N}\cap\Omega_{2,N}\cap\Omega_{3,N}$ we have that
\begin{align*}
&\forall j\in\mathcal{J}\colon \left\Vert H_N(t_j)\right\Vert \leqslant 2+\frac{(\log N)^{2}}{N^{1/4}},\\
&\forall j\in\mathcal{J}\colon \abs{\sprod{\vk e,H^k_N(t_j)\,\vk e}-\E{\sprod{\vk e,H^k_N(t_j)\,\vk e}}}< \mathfrak{C}^{k}\frac{(\log N)^{2k}}{\sqrt{N}}\\
&\forall j\in\mathcal{J}, t\in[t_j,t_{j+1}]\colon H_N(t)=H_N(t_j),\quad\text{or}\quad H_N(t) = H_N(t_{j+1}) ,
\end{align*}
implying that, on $\Omega_N$,
\begin{align*}
&\forall t\in[0,T]\colon \left\Vert H_N(t)\right\Vert \leqslant 2+\frac{(\log N)^{2}}{N^{1/4}},\\
&\forall t \in[0,T]\colon \abs{\sprod{\vk e,H_N^k(t)\,\vk e}-\E{\sprod{\vk e,H_N^k(t)\,\vk e}}}< \mathfrak{C}^{k}\frac{(\log N)^{2k}}{\sqrt{N}}.
\end{align*}
Finally, by combining \eqref{O1_bound}, \eqref{O2_bound} and \eqref{O3_bound},
\begin{align*}
\pk{\Omega_{N}^{\rm c}}\leqslant \pk{\Omega_{1,N}^{\rm c}} + \pk{\Omega_{2,N}^{\rm c}} + \pk{\Omega_{3,N}^{\rm c}}\leqslant 2e ^{-\frac{\nu^{\star}}{2}(\log N)^2} + e ^{-\frac{\nu^{\star}}{4}(\log N)^2} 
\leqslant e^{-\frac{\nu^{\star}}{8}(\log N)^{2}}.
\end{align*}
Hence, the claim follows by setting $\nu := \nu^{\star}/8$.
\QED

\medskip

We conclude this appendix by providing the proof of Lemma \ref{mu_full_expression}.

\medskip

\prooflem{mu_full_expression} Using \eqref{mu_recurrent_expansion}, as appearing in \nelem{lem_mu_recurrent_expansion}, we obtain on $\Omega_N$ that
\begin{align*}
    \mu^{\star}_{N}(t,N)&=\sqrt{Nq(t)}\sum_{k=0}^\infty\sprod{\vk e,\left(\frac{H_N(t)}{\mu^{\star}_{N}(t)}\right)^k\vk e}\\
    &=\sqrt{Nq(t)} + \sqrt{Nq(t)}\frac{\sprod{\vk e,H_N(t)\,\vk e}}{\mu^{\star}_{N}(t)} +\sqrt{Nq(t)}\sum_{k=2}^{\infty}\frac{\sprod{\vk e, H^k_N(t) \vk e}}{\bigl(\mu^{\star}_{N}(t)\bigr)^{k}}\\
    &=\sqrt{Nq(t)} + Z_{1,N}(t) + Z_{2,N}(t) =\sqrt{Nq(t)} + \tilde{\mathcal{G}}_N(t)\frac{(\log N)^2}{\sqrt{N}},
\end{align*}
where
\begin{align*}
    Z_{1,N}(t) &:= \sqrt{Nq(t)}
    \frac{\sprod{\vk e,H_N(t)\,\vk e}}{\mu^{\star}_{N}(t)},\\
    Z_{2,N}(t) &:= \sqrt{Nq(t)}\sum_{k=2}^{\infty}\frac{\sprod{\vk e, H^k_N(t)\,\vk e}}{\bigl(\mu^\star _{N}(t)\bigr)^{k}},\\
    \tilde{\mathcal{G}}_N(t) &:= \frac{\sqrt{N}}{(\log N)^2}(Z_{1,N}(t) + Z_{2,N}(t)).
\end{align*}
To find a uniform bound for $\tilde{\mathcal{G}}_N(t)$, we need first to find some 
uniform lower bound for $\mu_N^{\star}(t)$. 

Applying similar ideas as in the proof of \nelem{lem_mu_recurrent_expansion}, we can conclude that for sufficiently large $N$, we have that $\mu_N^{\star} > 7$. Utilizing this bound along with \eqref{A1} as given in \nelem{uniform_event}, and considering \eqref{mu_recurrent_expansion}, we find the following lower bound on $\Omega_N$:
\begin{align}
    \mu_N^{\star}(t) &= \sqrt{Nq(t)}\sum_{k=0}^\infty\sprod{\vk e,\left(\frac{H_N(t)}{\mu^{\star}_{N}(t)}\right)^k\vk e} \geqslant \sqrt{Nq(t)} -  \sqrt{Nq(t)}\sum_{k=1}^\infty\frac{\abs{\sprod{\vk e,H^k_N(t)\,\vk e}}}{\bigl(\mu^\star _{N}(t)\bigr)^{k}}\notag\\
    &\geqslant \sqrt{Nq(t)} -  \sqrt{Nq(t)}\sum_{k=1}^\infty\frac{3^k}{7^k}\geqslant \frac{\sqrt{Nq^-}}{4} \label{mu_rough_bound*}
\end{align}
Combining \eqref{mu_rough_bound*} with \eqref{A1},  we obtain that, for sufficiently large $N$ and for all $t\in[0,T]$, on $\Omega_N$,
\begin{align}
    \frac{\lvert\lvert H_N(t)\rvert\rvert}{\mu_N^{\star}(t)} \leqslant \frac{12}{\sqrt{Nq^-}}.\label{frac_rough_bound*}
\end{align}

We can now proceed with the bounds for the functions $Z_{1,N}(\cdot)$ and $Z_{2,N}(\cdot)$. Regarding the first function, using \eqref{mu_rough_bound*} and \eqref{A2} from \nelem{uniform_event} for $k=1$,
\begin{align}
\abs{Z_{1,N}(t)}\leqslant 4\mathfrak{C}\frac{(\log N)^{2}}{\sqrt{Nq^-}}\label{Z1_bound}
\end{align}
For $Z_{2,N}(t)$, relying on \eqref{frac_rough_bound*}, we have that for sufficiently large $N$,
\begin{align}
\abs{Z_{2,N}(t)}&\leqslant \sqrt{Nq^+}\sum_{k=2}^{\infty}\left(\frac{12}{\sqrt{Nq^-}}\right)^k\leqslant \frac{144\sqrt{q^+}}{q^-\sqrt{N}} \sum_{k=0}^{\infty}\left(\frac{12}{\sqrt{Nq^-}}\right)^k\leqslant  \frac{288\sqrt{q^+}}{q^-\sqrt{N}}\label{Z6_bound}
\end{align}
Combining \eqref{Z1_bound} and \eqref{Z6_bound}, we obtain that for sufficiently large $N$, 
\begin{align*}
\abs{\tilde{\mathcal{G}}_N(t)}&\leqslant \frac{4\mathfrak{C}}{\sqrt{q^-}} + \frac{288\sqrt{q^+}}{q^-}.
\end{align*}
Hence the claim follows.
\QED

\end{document}